\documentclass[dvips,aos,preprint]{imsart}

\RequirePackage[OT1]{fontenc}
\RequirePackage{amsthm,amsmath,amssymb,amsfonts}
\RequirePackage{latexsym}
\RequirePackage[numbers]{natbib}

\newtheorem{proposition}{Proposition}[section]
\newtheorem{theorem}{Theorem}[section]

\newtheorem{corollary}{Corollary}[section]
\newtheorem{lemma}{Lemma}[section]
\newtheorem{remark}{Remark}[section]

\arxiv{math.PR/0000000}

\startlocaldefs
\numberwithin{equation}{section}
\theoremstyle{plain}

\endlocaldefs

\begin{document}

\begin{frontmatter}
\title{On the Discrete Cram\'er-von Mises Statistics under Random Censorship}
\runtitle{On the Discrete Cram\'er-von Mises Statistics}

\begin{aug}
\author{\fnms{Dorival Le\~ao} \snm{}\thanksref{}\ead[label=e1]{leao@icmc.usp.br}}
\and
\author{\fnms{Alberto Ohashi} \thanksref{t3}
\ead[label=e3]{albertomfo@insper.edu.br}
}

\thankstext{t3}{Supported by CNPq Grant 798373}

\affiliation{Universidade de S\~ao Paulo and Insper Institute}

\address{Departamento de Matem\'atica Aplicada e Estat\'istica. \\
Universidade de S\~ao
Paulo, 13560-970\\
S\~ao Carlos - SP, Brazil.\\
\printead{e1}\\
}

\address{Insper Institute\\
04546-042  S\~ao Paulo SP, Brazil.\\
\printead{e3}\\
}
\end{aug}

\begin{abstract}
In this work, nonparametric log-rank-type statistical tests are introduced in order to verify homogeneity of purely discrete variables subject to arbitrary right-censoring for infinitely many categories. In particular, the Cram\'er-von Mises test statistics for discrete models under censoring is established. In order to introduce the test, we develop the weighted log-rank statistics in a general multivariate discrete setup which complements previous fundamental results of Gill~\cite{Gill} and Andersen et al.~\cite{ABGK}. Due to the presence of persistent jumps over the unbounded set of categories, the asymptotic distribution of the test is not distribution-free. The statistical test for a large class of weighted processes is described as a weighted series of independent chi-squared variables whose weights can be consistently estimated and the associated limiting covariance operator can be infinite-dimensional. The test is consistent to any alternative hypothesis and, in particular, it allows us to deal with crossing hazard functions. We also provide a simulation study in order to illustrate the theoretical results.
\end{abstract}

\begin{keyword}[class=AMS]
\kwd[Primary ]{62G10}
\kwd{62N03}
\kwd[; secondary ]{62N99 }
\kwd{62M99}
\end{keyword}

\begin{keyword}
\kwd{Nonparametric methods}
\kwd{Cram\'er-von Mises statistics}
\kwd{Survival analysis}
\kwd{Central Limit Theorem.}
\end{keyword}

\end{frontmatter}

\tableofcontents

\section{Introduction}

Discrete data analysis has great importance in several fields such as economics, biology, medicine, etc. Typically, discrete time data may occur either because the underlying data generating process is intrinsically discrete one or because they are discretely recorded. In most situations, statistical procedures based on continuous or mixed distributions cannot be directly applied to the purely discrete cases due to the lack of a continuum amount of information needed to validate the asymptotic results.

Frequently, one typically faces the problem of performing a data analysis along a time horizon subject to censoring, i.e., when some data at hand have occurred only within certain periods of time. To be more specifically, a general right-censoring scheme can be described as follows. Let $W_1^p , \ldots , W^p_{n_p}$ be independent positive random variables representing the survival times or times to some events of $n_p$ items in a population $p=1, \ldots , J$. The corresponding distribution function and intensity function of the $p$-th population are denoted by $F^p$ and $h^p$, respectively. A typical situation occurs when $\{W^p_m\}_{m=1}^{n_p}$ are censored on the right by independent positive random variables $\{C^p_m\}_{m=1}^{n_p}$. These censored variables $C^p_m$ are also assumed independent w.r.t $W^p_m$. Thus, in this general random censorship model one can only observe

\[
X_m^p = \min\{W^p_m, C^p_m\} , \quad \delta^p_m = 1\!\!1_{ \{ X^p_m = W_{m}^p \}}
\] where $\delta^p_m$ indicates whether $W^p_m$ is censored or not.

In a large number of applications it is of major importance to test the homogeneity of populations in the presence of censoring. A major problem in the literature is the development of nonparametric methods to test the null hypothesis
$$H_0 : F^1(\ell) =\cdots = F^J(\ell);~\ell\in \mathcal{X},$$
where $\mathcal{X}$ is the domain which encodes e.g the lifetime in a typical problem in survival analysis. The test should be consistent to a large class of alternative hypotheses. A lot of different test procedures have been proposed and studied so far~[see e.g~\cite{kalb} and other references therein]. One of the most important nonparametric methods to verify $H_0$ under random censorship is the well-known weighted log-rank test successfully developed in the setup of the Martingale theory proposed by Aalen~\cite{aalen} in survival analysis. See e.g~the works~[\cite{AlGBK},~\cite{ABG},~\cite{AET},~\cite{KM}~\cite{Gill},~\cite{FH}] and other references therein for this approach. The weighted log-rank test is one of the pillars of modern survival analysis. In fact, it is the most commonly used nonparametric test to compare two or more continuous populations with data that are subject to censoring. In spite of the large flexibility of the log-rank test, it may not exhibit good power to deal with non-proportional intensity functions. This fact is a major drawback of the log-rank theory and it is discussed in several applications with continuous distributions. See e.g~Klein and Moeschberger~\cite{KM} for a detailed discussion on this matter.

While there is a large number of works about nonparametric methods for continuous distributions subject to censoring~[see e.g~~\cite{FH},~\cite{AET} and~\cite{ABG} and other references therein], it is rather surprising that only few works investigate methods for lifetime discrete data. See e.g. the works~\cite{Gupta1},~\cite{GRI} and \cite{DK}. The study of nonparametric tests for $H_0$ and consistent w.r.t any alternative hypotheses for purely discrete distributions with full support under censoring is clearly rather important in many applications. This is the program we start to carry out in this work. We tackle this problem by naturally considering the weighted linear log-rank statistics given by

$$
LR_{q}({n^{\star}}, r) :=\sum_{\ell=1}^{r} \sum_{q_1\neq q}  U^{n_q}_{n_{q_1}} (n^{\star}, \ell ) \left[ \hat{h}^{n_q} (\ell) -  \hat{h}^{n_{q_1}} (\ell)\right],~n^\star\in \mathbb{N}^J,~r\ge 1,
$$ where $U^{n_q}_{n_{q_1}}$ is a suitable empirical weighted process, $\hat{h}^{n_p}$ is the Kaplan-Meier estimator for the intensity function $h^p$ and the variable $r$ encodes the number of categories which lives in an unbounded subset of the set of natural numbers $\mathbb{N}$. A large amount of attention in this article is devoted to the study of the asymptotic behavior of the following statistics

\begin{equation}\label{cvmintroduction}
CVM(n^{\star}) := \sum_{r=1}^{\infty} \sum_{q=1}^{J-1} LR_{q}^2({n^{\star}}, r) \hat{\phi}_q^2 (n^{\star}, r), \quad  n^{\star} \in \Bbb{N}^J,
\end{equation}
for suitable choices of weights $ \hat{\phi}_q$ and $U^{n_q}_{n_{q_1}}$. The statistics~(\ref{cvmintroduction}) is reminiscent from the classical Cram\'er-von Mises statistics largely used to compare continuous distributions under censoring. In this case, standard arguments based on Gaussian processes may be used to prove the correspondent asymptotic limit. See e.g ~\cite{Schumacher},~\cite{AET} and~\cite{KM} and other references therein. 
In the purely discrete case, the study of the asymptotic behavior of~(\ref{cvmintroduction}) is not trivial. In this case, there is no obvious (if any) underlying Gaussian process which describes the asymptotic behavior of~(\ref{cvmintroduction}). The reason is the occurrence of infinitely many persistent jumps as the sample size goes to infinity which causes some difficulties in establishing Lindeberg-type conditions.

Very general and fundamental asymptotic results for the log-rank statistics was established in the eighties by Gill~\cite{Gill} and Andersen et al~\cite{ABGK}. The key point of their asymptotic argument is based on the assumption that the size of the jumps of the underlying martingales goes to zero due to a Lindeberg-type condition. This of course does not hold when we are in a general discrete setup. In order to recover the discrete case, one needs to accommodate the relevant information on each category of the variable of interest. See e.g Murphy~\cite{murphy} and Lipster and Shiryaev~\cite{lipster} for a detailed discussion on this matter. To our best knowledge only few works have tried to tackle this problem. Gill \cite{Gill} proposed to spread the jumps at each category of the variable of interest over a neighborhood. His argument provides amenable time-changed processes whose their persistent jumps could in principle be asymptotically controlled under suitable conditions. See remarks after Theorem 4.2.1 in Gill~\cite{Gill}. However, his approach is far from being simple and it is not clear how one would apply his argument to a concrete discrete case in applications.

Along the lines developed by Gill~\cite{Gill} and Andersen et al.~\cite{ABGK}, other works~[see e.g~\cite{jones} and other references therein] study asymptotic properties of different classes of tests for $H_0$ but always restricted to continuous distributions. Another related work is Murphy~\cite{murphy} who proved a central limit theorem in a rather general discrete setup but the resulting weighted processes do not recover the log-rank statistics. Stute~\cite{stute} and Akritas~\cite{akritas} provide central limit theorems under censorship with distributions with atoms. They essentially study limit theorems for integrals w.r.t Kaplan-Meier estimators but they are not able to recover the log-rank setup.

In view of the remarks raised above, we derive one asymptotic result in the setup of the discrete log-rank theory towards applications to the statistics~(\ref{cvmintroduction}). Our approach is rather different from Gill~\cite{Gill}. We introduce suitable underlying discrete-time filtrations (over the samples) which allow us to get rid the unavoidable persistent asymptotic jumps of the underlying martingales. See Remark~\ref{fundremark} for more details. With these filtrations at hand, our first result (Theorem \ref{Conv_LR_infinity}) applied to the log-rank statistics is solely based on \textit{soft} arguments and well-known machinery from Martingale theory.

With the weighted log-rank asymptotic result at hand, we have prepared the basis for the asymptotic distribution (Theorem~(\ref{main_CVM})) of the Cram\'er-von Mises statistics in~(\ref{cvmintroduction}) for independent discrete populations under arbitrary right-censoring and with infinitely many categories. In particular, we provide a statistical test for $H_0$ which is consistent to any alternative hypothesis. Our test recovers previous works on discrete versions of Cram\'er-von Mises~\cite{LSS,chou,spinelli} for goodness of fit, without censoring and finite number of categories when applied to suitable choices of weighted processes. Moreover, the proposed test statistics supports a large class of empirical weighted processes which includes Fleming-Harrington~\cite{FH} and Tarone-Ware~\cite{TW} weights.

Contrary to the classical continuous case (see e.g~\cite{Schumacher}), the asymptotic distribution of the Cram\'er-von Mises statistic is not distribution-free because of the inherent nature of the jumps. Nevertheless, the underlying discrete structure allows us to write it as a weighted series of independent chi-squared variables with one degree of freedom where the weights can be consistently estimated from the data even when the limiting covariance operator is infinite-dimensional.

In order to illustrate the importance of our results to data analysis, we briefly compare classical nonparametric tests based on continuous distributions with the Cram\'er-von Mises test of this paper in a discretely recorded data set presenting crossing hazard functions. We show that the methodology developed in this paper allows us to detect crossing hazard functions when classical methodologies based on continuous distributions fail.


In order to deal with infinitely many categories, the limit theorems (see Theorems~\ref{main_asymptotic_result} and~\ref{main_CVM}) which describe the asymptotic behavior of~(\ref{cvmintroduction}) has to be worked out in the Hilbert space of square summable real sequences. In particular, Theorems~\ref{main_asymptotic_result} and~\ref{main_CVM} can be related to the work of Dedecker and Merlev\`ed \cite{DEDECKER}. In this work, they give necessary and sufficient conditions for a general stationary sequence of Hilbert space valued-random variables to satisfy the conditional central limit theorem. In particular they apply their results to characterize the asymptotic distribution of the Cram\'er-von Mises statistics essentially involving the empirical distribution function for continuous and/or discrete random variables. Unfortunately, their results do not recover censored random variables. Our central limit theorem complements~\cite{DEDECKER} in the particular case of discrete variables arising in a right-censorship model where the empirical distribution function is not available for analysis.

The remainder of this paper is organized in the following way. After fixing the notation and recalling some basics of the inference of the multiplicative intensity model in Section~\ref{multmodel}, we formulate and state the main asymptotic results of this article in Sections~\ref{firstsection},~\ref{log-ranksection} and~\ref{CVMsection}. Section~\ref{simsection} reports a simulation study and the proofs of Theorems~\ref{Conv_LR_infinity} and~\ref{main_asymptotic_result}, Proposition~\ref{HAass} and Theorem~\ref{main_CVM} are reported in Sections~\ref{proofTH1},~\ref{proofTH2},~\ref{proofProp} and~\ref{proofTH3}, respectively. The Appendix contains some technical results used in Theorem~\ref{Conv_LR_infinity}.

\section{Inference for the Multiplicative Intensity Model for Discrete Random Variables}\label{multmodel}
In this section, we introduce the basic notation and the discrete model used in this paper. Let $\mathbb{N}$ be the set of non-negative integers and let $W$ and $C$ be two independent discrete $\mathbb{N}$-valued random variables. The discrete random variable $W$ describes the event of interest, while the discrete random variable $C$ denotes the censoring variable. The reader may think $C$ as the random variable which describes the censoring in a given statistical problem. Let us fix an underlying probability space $(\Omega,\mathcal{F},\mathbb{P})$. Throughout this paper, we make use of the following notation: if $Z$ and $Y$ are two adapted processes, we write $(Y.Z)(i):= Y(0)Z(0) + \sum_{\ell=1}^i Y(\ell)\Delta Z(\ell)$ for $i\ge 1$, where $\Delta Y(\ell)= Y(\ell)-Y(\ell-1)$. In order to distinguish discrete-time stochastic integrals from simple Riemman sums, if $N$ is a predictable process and $M$ is a martingale, we shall write the correspondent discrete-time stochastic integral as $\int_0^i N(\ell)dM(\ell)= (N.M)(i)$. Moreover, if $M$ is discrete-time martingale then $\langle M\rangle$ denotes the usual predictable bracket.

We fix once and for all, a natural number $J\geq 1$ and we write $\mathcal{J} = \{1,\ldots, J\}.$ In the remainder of this paper, we denote $n^\star:=(n_1,\ldots,n_J)\in\mathbb{N}^J$ and $N_J:=\max\{n_1,\ldots, n_J\}$. Here $n_q$ denotes the size of the random sample at hand correspondent to the population $q\in \mathcal{J}$. The symbols $\vee$ and $\wedge$ will denote the maximum and minimum between real numbers, respectively. Let us now describe the general right-censorship discrete random model of this work. For a given $p\in \mathcal{J}$, let $(W^p,C^p)$ be a given population where $W^p$ and $C^p$ are independent discrete random variables which must be interpreted as $W$ and $C$, respectively. Let $X^p$ be a discrete random variable defined by

$$
 X^p := W^p \wedge C^p,\quad p\in \mathcal{J}.
$$
For any $p \in \mathcal{J}$, we take independent random samples $\{ (W_1^p , C_1^p) , \cdots , (W_{n_p}^p , C_{n_p}^p) \}$ from the population $(W^p ,C^p)$ where $1\le n_p < \infty$. To shorten notation, throughout this paper we assume that $\mathbb{P}[C^p=0]=\mathbb{P}[W^p=0]=0$ for any $p\in \mathcal{J}$. With these random samples at hand, we introduce

$$ X_m^p := W_m^p \wedge C_m^p, \quad V^{p}_m(i):=1\!\!1_{\{X_{m}^{p} \ge i\}}$$

\noindent and the counting processes

\[
R^p_m(i) := 1\!\!1_{\{X_m^p \leq i, X_m^p=W_m^p\}} \quad  \text{and} \quad R^{C,p}_m(i) := 1\!\!1_{\{X_m^p \leq i, X_m^p=C_m^p\}},
\]

\noindent for  $ m =1 , \ldots , n_p$ and $i\ge 0$. The counting processes associated with the $p$-th random sample are given by

\begin{equation}\label{doobsample}
R^{n_p}(i) := \sum_{m=1}^{n_p} R^p_m(i) \quad  \text{and} \quad R^{C,n_p}(i) := \sum_{m=1}^{n_p} R^{C,p}_m(i), \quad i\ge 0.
\end{equation}

In order to take into account all the information generated by the random sample at hand, we define

$$\mathcal{F}^{n_p}_i:=\bigvee _{m=1}^{n_p}\mathcal{A}^p_{m,i},$$

\noindent where

$$\mathcal{A}_{m,i}^p := \sigma(\Delta R^{C,p}_m(\ell), \Delta R^p_m(\ell);1\le \ell\le i),$$

\noindent for each $i\ge 1$ and $ m =1 , \ldots , n_p$. We also set $\mathcal{F}^{n_p}_0:=\mathcal{A}^p_{m,0}:=\{\emptyset, \Omega\}$ for $m=1,\ldots, n_p$ and $p\in \mathcal{J}$. Here, for a given family of random variables $D$ the class $\sigma(D)$ is the smallest sigma-algebra making all of $D$ measurable. Moreover, $\bigvee$ is the smallest sigma-algebra generated by a union of sigma-algebras.

In this paper, several types of filtrations will play different rules. The filtration

$$\mathbb{F}^{n_p}:=\{\mathcal{F}^{n_p}_i;i\ge 0\}$$
will be used to perform the Doob-Meyer decomposition for the counting process $R^{n_p}$ as follows. By the very definition, for each $m\in \{1,\ldots, n_p \}$ and $p\in \mathcal{J}$ the Doob-Meyer decomposition of $R^p_m$ w.r.t $\mathbb{A}^{p}_m:= \{\mathcal{A}^{p}_{m,i} , i\ge 0\}$ reads

$$R^p_m(i)=Y^p_m(i) + N^p_m(i),~i\ge 1$$
where $N^p_m$ is the compensator of $R^p_m$ and $Y^p_m$ is the associated martingale. Here $R^p_m(0)=Y^p_m(0)=N^p_m(0)=0$ a.s. The definition of the filtration $\mathbb{F}^{n_p}$ yields the following Doob-Meyer decomposition

$$
R^{n_p}(i) = \sum_{m=1}^{n_p} \Big(Y^{p}_m(i) + N^{n_p}_m(i)\Big) =: Y^{n_p}(i) + N^{n_p}(i),~i\ge 1.
$$

\noindent By the very definition, the following $\mathbb{F}^{n_p}$-decomposition holds

\begin{equation} \label{doobsample}
R^{n_p}(i) = Y^{n_p}(i) + \sum_{\ell=1}^{i} V^{n_p} (\ell) \Delta H^p (\ell)  \quad i\ge 1,
\end{equation}

\noindent where

$$V^{n_p}(\cdot) := \sum_{m=1}^{n_p} V^{p}_m(\cdot) \quad \text{and} \quad  H^p(\cdot):= \sum_{j=1}^{\cdot} h^p(j).$$

\noindent Here $h^p$ is the intensity function associated to the random variable of interest $W^p$ which can be written

\begin{equation}\label{hazard}
h^p (j) := \frac{\mathbb{P} [ W^p = j ]}{\mathbb{P} [ W^p \ge j ]},\quad j\ge 0,
\end{equation}
provided that $\mathbb{P}[W^p\ge j]>0$.
\begin{remark}\label{covarYrem}
By construction, one should notice that each $\mathbb{A}^{p}_m$-martingale $Y^p_m$ is also a martingale w.r.t $\mathbb{F}^{n_p}$ for every $ m =1 , \ldots , n_p$. Moreover, the martingales $Y^p_m$ and $Y^p_j$ are independent for any $m \neq j$ in $\{1 , \ldots , n_p\}$. As a consequence, the following representation for the predictable bracket holds

$$
\langle Y^{n_p} \rangle (i) =  \sum_{\ell=1}^{i} V^{n_p}(\ell) h^p(\ell) \left[ 1-h^p(\ell)\right],\quad i\ge 1.
$$
\end{remark}
\noindent In the remainder of this paper, we denote

$$\theta^p(\ell):=\mathbb{P}[X^p\ge \ell];~p\in \mathcal{J},\ell\ge 0,$$
and whenever necessary, we can always assume that for each $q\in \mathcal{J}$ there exists a category $i$ such that $\theta^q(i)>0$. This ensures for instance that we can write $h^q$ as in~(\ref{hazard}) on $\{1,\ldots, i\}$. One can easily check the following elementary property.

\begin{remark}\label{binomialrem}
$V^{n_p}(\ell)$  has binomial distribution with parameters $n_p$ and $\theta^p(\ell)$ for each $\ell\ge 1$. The conditional distribution of $\Delta R^{n_p}(\ell)$ given $V^{n_p}(\ell) = j$ is binomial with parameters $j$ and $h^p(\ell)$, for every $\ell\ge 1$ and $j=1, \ldots , n_p$. Moreover,

\begin{equation}\label{intensitysample}
\mathbb{E}[\Delta R^{n_p}(i) \mid V^{n_p}(i) ] ~ = ~ V^{n_p}(i) h^p(i) = \mathbb{E}[\Delta R^{n_p}(i) \mid \mathcal{F}^{n_p}_{i-1} ], \quad i\ge 1.
\end{equation}
\end{remark}

We are now in position to derive an empirical estimator for the intensity function $h^p$. The idea is fully based on the relations~(\ref{doobsample}) and~(\ref{intensitysample}). The martingale component $Y^{n_p}$ is interpreted as a noise and the predictable component $N^{n_p}$ contains all the information about the law of the discrete random variable $W^p$ needed for the estimation. Therefore, it is natural to introduce an estimator for $h^p$ based on the condition that $Y^{n_p}=0$. In fact, if $Y^{n_p}=0$, we recover the Kaplan-Meier estimator as follows

\begin{equation}\label{kaplanmeyer}
\hat{h}^{n_p}(i) = \frac{\Delta R^{n_p}(i)}{V^{n_p}(i)}1\!\!1_{\{V^{n_p}(i)>0 \}} ;\quad i\ge 1.
\end{equation}




Based on~(\ref{kaplanmeyer}), we then define the following estimators for $H^p$ and for the law $\pi^p$ of the discrete random variable $W^p$, respectively, as follows

\begin{equation}\label{eq_estimador_H}
\hat{H}^{n_p}(i):= \sum_{\ell=1}^{i}\hat{h}^{n_p}(\ell),
\end{equation}

\begin{equation}\label{eq estimator W}
\hat{\pi}_i^{n_p}:= \hat{h}^{n_p}(i)\prod_{\ell=1}^{i-1}[1-\hat{h}^{n_p}(\ell)].
\end{equation}

\noindent for $i\ge 1$. The Doob-Meyer decomposition~(\ref{doobsample}) and~(\ref{kaplanmeyer}) yield

\begin{equation}\label{eq fundam}
\hat{H}^{n_p} (i) - H^{p} (i) = \int_{0}^i \frac{1}{V^{n_p}(\ell)}1\!\!1_{\{V^{n_p}(\ell)>0 \}} d Y^{n_p}(\ell) + r_{n_p} (i); i\ge 1,
\end{equation}

\noindent where the predictable component $r_{n_p}(i):=\sum_{\ell=1}^{i}h^p(\ell)[1\!\!1_{\{V^{n_p}(\ell)>0\}}-1]$ vanishes as $n_p\rightarrow \infty$ due to a Borel-Cantelli argument provided $\theta^p(i) > 0$. The following useful remark gives the asymptotics of the above estimators. These results can be easily proved by routine arguments so we omit the details.

\begin{remark}\label{hestimator}
 For a given $p\in \mathcal{J}$, let $i$ be a positive integer such that $\theta^p(i)>0$. From identity~(\ref{eq fundam}) we may conclude that $\hat{H}^{n_p}\rightarrow H^p$ uniformly (over $\{1,\ldots,i \}$) in probability as $n_p\rightarrow \infty$. Moreover, $\hat{h}^{n_p}(\ell) \rightarrow h^p(\ell)$ and $\hat{\pi}_\ell^{n_p} \rightarrow \pi_\ell^p$ in probability as $n_p\rightarrow \infty$, for each $\ell\in \{1,\ldots,i \}$. If $\theta^p(i)>0$ for every $i\ge 1$, then $\sup_{\ell\ge 1}|\hat{h}^{n_p}(\ell)-h^p(\ell)|\rightarrow 0$ in probability as $n_p\rightarrow \infty$.


\end{remark}

\section{Asymptotic Distribution for Discrete Stochastic Integrals}\label{firstsection}
In this section, we provide the asymptotic results which will be the basis for the statistical tests in this article. From~(\ref{eq fundam}) and Remark~\ref{hestimator}, we know that the accumulated intensity process $H^{p}$ admits a natural class of consistent estimators $\hat{H}^{n_p}$ in such way that $\hat{H}^{n_p}- H^{p}$ is a discrete-time $\mathbb{F}^{n_p}$-semimartingale of the form~(\ref{eq fundam}). Therefore, a natural strategy will be based on a martingale central limit theorem. In the previous section, we have defined the filtration family $\{ \mathbb{F}^{n_p}; p \in \mathcal{J},n_p\ge 1\}$ where the Doob-Meyer decomposition and the resulting Kaplan-Meier estimator were performed. In the sequel, in order to obtain the asymptotic distribution, we are forced to use different types of filtrations. Let

$$R^{n^{\star}}(i):= \sum_{k=1}^JR^{n_k}(i);~i\ge 0,$$

\noindent be the total number of events of interest  at category $i$ and let

$$V^{n^{\star}}(i) :=\sum_{k=1}^J V^{n_k}(i);~i\ge 0,$$

\noindent be the total number at risk at category $i$. In order to keep track the limiting martingale behavior at different samples, we introduce the filtration $\mathbb{F}=\{\mathcal{F}_i;i\ge 0\}$ generated by the whole information available at each category as follows

$$\mathcal{F}_i:=\bigvee_{n_p;~p\in \mathcal{J}}\mathcal{F}^{n_p}_i\quad i\ge 0.$$

\noindent In the sequel, to shorten notation we introduce

$$\mathcal{V}^p_m(j) := \big(V^p_1(j),\ldots, V^p_m(j)\big), \quad \mathcal{R}^p_m(j) := \big(R^p_1(j), \ldots, R^p_m(j)\big),$$

\noindent and

$$\mathcal{R}^{C,p}_m(j) := \big(R^{C,p}_1(j), \ldots, R^{C,p}_m(j)\big),\quad 1\le m \le N_J; j\ge 1,~p\in\mathcal{J}.$$

Let us now introduce another filtration which will be the basis for our asymptotic results. This filtration family is carefully chosen as follows. For a given category $j\ge 1$ and $n^\star\in \mathbb{N}^J$, we define the filtration $\mathcal{G}^{n^\star}(j):=\{\mathcal{G}^{n^\star}_m(j); 0\le m\le N_J\}$ along the samples as follows

\begin{eqnarray*}
\mathcal{G}^{n^\star}_m(j)&:=&\sigma \{\big(V^{n_p}(j), \Delta R^{n^\star}(j-1), \mathcal{V}^p_{m+1}(j),\Delta \mathcal{R}^p_m(j),\Delta \mathcal{R}^{C,p}_m(j)\big),~ p \in \mathcal{J} \},
\end{eqnarray*}
for any $0\le m\le N_J-1$ and

\begin{eqnarray*}
\mathcal{G}^{n^\star}_{N_J}(j)&:=&\sigma \{ \big(V^{n_p}(j), \Delta R^{n^\star}(j-1), \mathcal{V}^p_{N_J}(j),\Delta \mathcal{R}^p_{N_J}(j),\Delta \mathcal{R}^{C,p}_{N_J}(j)\big), ~ p \in \mathcal{J}\},
\end{eqnarray*} Here we set $\Delta R^{n^\star}(0)=R^{n^\star}(0)$ and $\Delta \mathcal{R}^p_0=\Delta \mathcal{R}^{C,p}_0=0$.
\begin{remark}\label{fundremark}
The building block for the asymptotic results of this article is based on a martingale structure over the filtration $\mathcal{G}^{n^\star}(j)$ along the samples for a given category $j\ge 1$. In fact, one can readily see that for every $p\in\mathcal{J}$, $n^\star\in\mathbb{N}^J$ and $j\ge 1$, $Y^p_\cdot(j)$ is a $\mathcal{G}^{n^\star}(j)$-martingale array difference. Moreover, for every $p\in \mathcal{J}$ and $n^\star\in\mathbb{N}^J$, $R^{n_p}(j)$ has the same Doob-Meyer decomposition w.r.t. $~\mathcal{G}^{n^\star}(j)$ (over the samples) and $\mathbb{F}^{n_p}$ (over the categories). In order to recover the classical Tarone-Ware and Harrington-Fleming weighted processes, we include the term $\Delta R^{n^\star}(j-1)$ in the definition of the filtration.
This structure allows us to accommodate the persistent jumps at each category as the sample size goes to infinity. This strategy is rather different from Gill~\cite{Gill} who used a time-changed argument on the level of categories in order to deal with the jumps of partially discrete distributions. Gill's idea is to spread the jump of $R^{n_p} (j)$ at category $j$ over a time interval which is inserted at this category.
\end{remark}

In the remainder of this paper, we always consider the random variables $1/V^{n_q}(\ell)$ as been multiplied by indicator functions $1\!\!1_{\{V^{n_q}(\ell)>0 \}}$ for any $q\in \mathcal{J}$ and $\ell\ge 1$.

\subsection{A Martingale Central Limit Theorem} Let us now describe a list of the technical assumptions which will constitute the basis for the asymptotic results of this section. For any pair $q_1 \neq q$ in $\mathcal{J}$ and $n^\star=(n_1,\ldots,n_J)\in\mathbb{N}^J$, we are going to write $$U_{n_{q_1}}^{n_q}(n^\star,\cdot)=\{U_{n_{q_1}}^{n_q}(n^\star,i); i\ge 1 \}$$
as an $\mathbb{F}$-predictable process satisfying some technical assumptions. In the remainder of this work, $n^\star\rightarrow \infty$ means $n_p\rightarrow \infty$ for every $p\in \mathcal{J}$.

\

\noindent \textbf{(M1)} For each $(n_q,n_{q_1})\in\mathbb{N}^2$ and $n^\star\in\mathbb{N}^J$, $U_{n_{q_1}}^{n_q}(n^\star,i)$ is $\mathcal{G}^{n^\star}_{0}(i)$-measurable for every $i\ge 1$;

\

\noindent \textbf{(M2)} There exists $\delta>0$ such that

$$\lim_{n^\star \rightarrow \infty} {n_q} \Bigg|\frac{U_{n_{q_1}}^{n_q}(n^\star,i)}{V^{n_q}(i)}\Bigg|^{2+\delta}~\vee~n_{q_1} \Bigg|\frac{U_{n_{q_1}}^{n_q}(n^\star,i)}{V^{n_{q_1}}(i)}\Bigg|^{2+\delta} =0$$

\noindent in probability for each $i\ge 1$;

\

\noindent \textbf{(M3)} For any $q_2 \in \{q,q_1\}$ and $\ell \ge 1$ there exists a constant $\alpha^{q_2}_{q,q_1}(\ell)$ such that

$$\frac{|U^{n_{q}}_{n_{q_1}}(n^\star,\ell)|^2}{V^{n_{q_2}}(\ell)}\rightarrow \alpha^{q_2}_{q,q_1}(\ell)$$

\noindent in probability as $n^\star\rightarrow \infty$.

\

\noindent \textbf{(M3$^{\prime}$)} For any $q_2 \in \{q,q_1\}$ we have that

$$ \sum_{\ell=1}^{\infty} \limsup_{n^\star}\mathbb{E}\frac{|U^{n_{q}}_{n_{q_1}}(n^\star,\ell)|^2}{V^{n_{q_2}}(\ell)}< \infty$$
and
$$
\sum_{\ell=1}^{\infty}\Bigg|\frac{|U^{n_{q}}_{n_{q_1}}(n^\star,\ell)|^2}{V^{n_{q_2}}(\ell)}-\alpha^{q_2}_{q,q_1}(\ell)\Bigg|\rightarrow 0$$
in probability as $n^\star\rightarrow \infty$.

\

\noindent \textbf{(M4)} For any $q_2\in \mathcal{J}$~$(q_2\neq q,q_2\neq q_1)$ and $\ell\ge 1$ there exists a constant $\beta^q_{q_1,q_2}(\ell)$ such that

$$\frac{U^{n_{q}}_{n_{q_1}}(n^\star,\ell)U^{n_q}_{n_{q_2}}(n^\star,\ell)}{V^{n_q}(\ell)}\rightarrow \beta^{q}_{q_1,q_2}(\ell)$$

\noindent in probability as $n^\star\rightarrow \infty$.

\

\noindent \textbf{(M4$^{\prime}$)} For any $q_2\in \mathcal{J}$~$(q_2\neq q,q_2\neq q_1)$ we have that

$$ \sum_{\ell=1}^{\infty} \limsup_{n^\star}\mathbb{E}\frac{\big|U^{n_{q}}_{n_{q_1}}(n^\star,\ell)U^{n_q}_{n_{q_2}}(n^\star,\ell)\big|}{V^{n_q}(\ell)} < \infty
$$
and
$$ \sum_{\ell=1}^{\infty}\Bigg|\frac{U^{n_{q}}_{n_{q_1}}(n^\star,\ell)
U^{n_q}_{n_{q_2}}(n^\star,\ell)}{V^{n_q}(\ell)} -\beta^q_{q_1,q_2}(\ell)
\Bigg|\rightarrow 0
$$
in probability as $n^\star\rightarrow \infty$.

When the number of categories is finite then \textbf{(M3$^{\prime}$)} and \textbf{(M4$^{\prime}$)} are not necessary. Since one of the main applications of the theoretical results of this paper lies in lifetime data analysis under censoring, then it is crucial to work under the setup of infinitely many categories. Assumptions~\textbf{(M3$^{\prime}$)} and~\textbf{(M4$^{\prime}$)} encode exactly this situation.

In fact, we are going to show that a large class of weighted processes satisfy the above list of technical assumptions. See Section~\ref{weighted} for more details. Let us now define a family of random variables which encodes any hypothesis test related to the homogeneity of several discrete populations under censoring with infinitely many categories. For a given $n^\star\in\mathbb{N}^J$ and $q\in \mathcal{J}$, we define the following  random variables

\begin{equation}\label{marraydif}
\xi^{n^{\star}}_{m,q}(\ell):=  \sum_{q_1 \neq q}   U_{n_{q_1}}^{n_q}(n^\star,\ell) \left[ \frac{\Delta Y^{q}_m(\ell)}
{V^{n_{q}}(\ell)} -  \frac{\Delta Y^{{q_1}}_m(\ell)}{V^{n_{q_1}}(\ell)} \right],
\end{equation}
where $m=1,\ldots, N_J$ and we set $\Delta Y^{k}_m=0$ if $n_k < m \le N_J$ and $k\in \mathcal{J}$.

A first simple remark is that $\{ \xi^{n^{\star}}_{m,q}(\ell); 1\le m \le N_J\}$ is a martingale-array difference w.r.t the filtration $\mathcal{G}^{n^\star} (\ell)$ (see Lemma \ref{multimardiff}) for each $q\in \mathcal{J}$, $n^\star\in \mathbb{N}^J$ and $\ell\ge 1$. Thus, we shall apply usual arguments from martingale theory (see Lemma \ref{con_ell_2}) to prove that the sequence $\sum_{m=1}^{N_J} \xi^{n^{\star}}_{m,q}(\ell)$ converges weakly to a zero mean Gaussian distribution with variance $\phi^2_q(\ell)$ given by

\begin{eqnarray}
\nonumber\phi^2_q(\ell)&:=&\sum_{q_1\neq q}\alpha_{q,q_1}^{q_1}(\ell)h^{q_1}(\ell)[1-h^{q_1}(\ell)] + \sum_{q_1\neq q}\alpha_{q,q_1}^{q}(\ell)h^{q}(\ell)[1-h^{q}(\ell)] \\
\nonumber & &\\
\label{asvar1} &+& 2\sum_{(q_1,q_2)\in A_q} \beta^q_{q_1,q_2}(\ell)h^q(\ell)[1-h^q(\ell)]; \quad \ell \ge 1,
\end{eqnarray}
where the family of functions $\alpha^{q_1}_{q,q_1},~\alpha^q_{q,q_1}, \beta^q_{q_1,q_2}:\mathbb{N}\rightarrow \mathbb{R}$ in~(\ref{asvar1}) are given in \textbf{(M3)} and \textbf{(M4)}, respectively. In~(\ref{asvar1}), we denote $A_q:=\{(x,y)\in \mathcal{J}\times \mathcal{J}; x\neq y, x\neq q, y\neq q, 1\le x < y \le J \}$ for $q\in \mathcal{J}$.

In Lemma~\ref{con_ell_2}, we also show that the asymptotic variance $\phi^2_q(\ell)$ can be consistently estimated by

\begin{eqnarray}
\nonumber\hat{\phi}^2_{q,n^\star}(\ell)&:=& \sum_{q_1\neq q}\Bigg[\frac{|U_{n_{q_1}}^{n_q}(n^\star,\ell)|^2}{V^{n_{q_1}}(\ell)}  \nonumber\hat{h}^{n_{q_1}} (\ell) [1-\hat{h}^{n_{q_1}} (\ell)] + \frac{|U_{n_{q_1}}^{n_q}(n^\star,\ell)|^2}{V^{n_q}(\ell)}  \nonumber\hat{h}^{n_{q}} (\ell) [1-\hat{h}^{n_{q}} (\ell)]\Bigg] \\
\label{varest}& &\\
\nonumber&+& 2\sum_{(q_1,q_2)\in A_q}\frac{U^{n_q}_{n_{q_1}}(n^\star,\ell)U^{n_q}_{n_{q_2}}(n^\star,\ell)}{V^{n_q}(\ell)}\hat{h}^{n_{q}} (\ell) [1-\hat{h}^{n_{q}} (\ell)],
\end{eqnarray}

\noindent for $n^\star\in\mathbb{N}^J$ and $\ell\ge 1$.

In the sequel, the analysis will be based on the following multi-dimensional process

$$\xi^{n^\star}(\ell) := \Bigg(\sum_{m=1}^{N_J} \xi^{n^\star}_{m,1}(\ell),\ldots, \sum_{m=1}^{N_J} \xi^{n^\star}_{m,J}(\ell)\Bigg);\quad \ell\ge 1.$$
The multi-dimensional case requires additional assumptions on a given weighted process. Given $r\neq k$ in $\mathcal{J}$, let $U$ be a weighted $\mathbb{F}$-predictable process which satisfies the following assumptions:

\

\noindent \textbf{(H1)} For any $q_1\neq k$ and $\ell\ge 1$, there exist constants $\gamma^{k,r}_{q_1}(\ell)$ and $\eta^{k,r}_{q_1}(\ell)$ such that

$$\frac{U^{n_k}_{n_{q_1}}(n^\star,\ell) U^{n_r}_{n_{q_1}}(n^\star,\ell)}{V^{n_{q_1}}(\ell)}\rightarrow \gamma^{k,r}_{q_1}(\ell)$$

$$
\frac{U^{n_k}_{n_{q_1}}(n^\star,\ell) U^{n_r}_{n_{k}}(n^\star,\ell)}{V^{n_{k}}(\ell)}\rightarrow \eta^{k,r}_{q_1}(\ell)
$$

\noindent in probability as $n^\star\rightarrow \infty$. Moreover, $U^{n_q}_{n_k}U^{n_q}_{n_r}$ is non-negative a.s for every $q\in\mathcal{J}$ with $q\ne k$ and $q\neq r$.

\

\noindent \textbf{(H1$^{\prime}$)} For any $q_1\neq k$, we assume that

$$\sum_{\ell=1}^{\infty} \limsup_{n^\star}\mathbb{E}\Bigg\{\frac{\big|U^{n_k}_{n_{q_1}}(n^\star,\ell) U^{n_r}_{n_{q_1}}(n^\star,\ell)\big|}{V^{n_{q_1}}(\ell)} + \frac{\big|U^{n_k}_{n_{q_1}}(n^\star,\ell) U^{n_r}_{n_{k}}(n^\star,\ell)\big|}{V^{n_{k}}(\ell)}\Bigg\} < \infty
$$
and

$$
\sum_{\ell=1}^{\infty}\Bigg|\frac{U^{n_k}_{n_{q_1}}(n^\star,\ell) U^{n_r}_{n_{k}}(n^\star,\ell)}{V^{n_{k}}(\ell)}-\eta^{k,r}_{q_1}(\ell)\Bigg|+
\sum_{\ell=1}^{\infty}\Bigg|\frac{U^{n_k}_{n_{q_1}}(n^\star,\ell) U^{n_r}_{n_{q_1}}(n^\star,\ell)}{V^{n_{q_1}}(\ell)}-\gamma^{k,r}_{q_1}(\ell)\Bigg|\rightarrow 0
$$
\noindent in probability as $n^\star\rightarrow \infty$.

Of course, the above assumptions \textbf{(H1-H1$^{\prime}$)} only make sense if $J\ge 3$. Without any loss of generality, throughout this section we assume that $J\ge 3$. See Remark~\ref{j<3}. In the sequel, if $k\neq r$ in $\mathcal{J}$ then we set $A(k,r):= \{q_1\in \mathcal{J}; q_1\neq k,q_1\neq r\}$, and we denote

\begin{eqnarray} \label{conv_operator}
\nonumber\psi(k,r,\ell)&:=&\sum_{q_1 \in A(k,r)}\gamma^{k,r}_{q_1}(\ell)h^{q_1}(\ell)[1-h^{q_1}(\ell)]-\sum_{q_1\neq k}\eta^{k,r}_{q_1}(\ell) h^{k}(\ell)[1-h^{k}(\ell)]\\
\nonumber & &\\
&-& \sum_{q_2\neq r}\eta^{r,k}_{q_2}(\ell)h^{r}(\ell)[1-h^{r}(\ell)], \quad \ell\ge 1,
\end{eqnarray}

\noindent where the functions $\gamma^{k,r}_{q_1}, \eta^{k,r}_{q_1},  \eta^{r,k}_{q_2}:\mathbb{N}\rightarrow \mathbb{R}$ are defined via a weighted process satisfying assumptions \textbf{(H1)}. Since $\gamma^{k,r}_{q_1}(\cdot) = \gamma^{r,k}_{q_1}(\cdot)$ for every $k\neq r$ in $\mathcal{J}$ provided $q_1\in A(k,r)$ then the symmetrization in the second line of~(\ref{conv_operator}) yields $\psi(k,r,\cdot) = \psi(r,k,\cdot)$. Moreover, for $k\neq r$ we denote

\begin{eqnarray}
\nonumber \hat{\psi}_{n^\star}(k,r,\ell)&:=&\sum_{q_1 \in A(k,r)}\frac{U^{n_k}_{n_{q_1}}(n^\star,\ell) U^{n_r}_{n_{q_1}}(n^\star,\ell)}{V^{n_{q_1}}(\ell)} \hat{h}^{q_1}(\ell)[1-\hat{h}^{q_1}(\ell)]\\
\nonumber & &\\
\nonumber &-&\sum_{q_1\neq k}\frac{U^{n_k}_{n_{q_1}}(n^\star,\ell) U^{n_r}_{n_{k}}(n^\star,\ell)}{V^{n_{k}}(\ell)} \hat{h}^{k}(\ell)[1-\hat{h}^{k}(\ell)]\\
\nonumber & &\\
\label{matriz1}&-& \sum_{q_2\neq r}\frac{U^{n_r}_{n_{q_2}}(n^\star,\ell) U^{n_k}_{n_{r}}(n^\star,\ell)}{V^{n_{r}}(\ell)} \hat{h}^r(\ell)[1-\hat{h}^{r}(\ell)],
\end{eqnarray}

\noindent for $n^\star\in\mathbb{N}^J,~\ell\ge 1$. Again symmetrization in the second line of~(\ref{matriz1}) yields $\hat{\psi}_{n^\star}(k,r,\cdot)=\hat{\psi}_{n^\star}(r,k,\cdot)$ a.s for every $k\neq r$ in $\mathcal{J}$ and $n^\star\in \mathbb{N}^J$. Now we are in position to state the first result of this section. In the sequel, we set $\Gamma(0)=0$ and we denote $\Gamma(i):= \sum_{\ell=1}^iQ(\ell),~i\ge 1$, where $Q$ is the self-adjoint operator defined by the following quadratic form

\begin{equation}\label{covoperator}
\langle Q(\ell)a,a\rangle_{\mathbb{R}^J} =\sum_{k=1}^J a^2_k\phi^2_k(\ell) + 2\sum_{1\le r< k \le J} a_ra_k \psi(k,r,\ell);~a\in \mathbb{R}^J,~\ell\ge 1.
\end{equation}
Convergence stated in~(\ref{c2})~and~(\ref{non2}) ensure that the quadratic form~(\ref{covoperator}) is actually non-negative. We also define the self-adjoint random operator $\hat{Q}(n^\star,\ell)$ induced by the quadratic form

\begin{equation}\label{covoperatorest}
\langle \hat{Q}(n^\star,\ell) a,a\rangle_{\mathbb{R}^J} =\sum_{k=1}^J a^2_k \hat{\phi}^2_{k,n^\star}(\ell) + 2\sum_{1\le r< k \le J} a_ra_k \hat{\psi}_{n^\star}(k,r,\ell);~a\in \mathbb{R}^J,~\ell\ge 1.
\end{equation}
We set $\hat{\Gamma}(n^\star, 0):=0$ and $\hat{\Gamma}(n^\star,i):=\sum_{\ell=1}^i\hat{Q}(n^\star,\ell);~n^\star\in\mathbb{N}^J,~i\ge 1$.

Since the variables of interest $\{W^p; p\in \mathcal{J}\}$ assume values in an unbounded set in a typical discrete lifetime data, it is important to introduce the following objects. Let $\{d^l_{n^*};n^\star\in\mathbb{N}^J\}$ and $\{d^u_{n^*};n^\star\in\mathbb{N}^J\}$ be two sequences of $\mathbb{F}$-stopping times which satisfy the following hypotheses:

\

\noindent \textbf{(S1)} $d^l_{n^*} < d^u_{n^*} < \infty$ a.s for every $n^*\in \mathbb{N}^J$ and there exists a pair $(d^l,d^u)\in \bar{\mathbb{N}}^2$ such that $1 \leq d^l < d^u \leq \infty$, and

$$d^l_{n^*} \rightarrow d^l~\text{and}~d^u_{n^*} \rightarrow d^u~\text{in probability}~\text{as}~n^\star\rightarrow \infty.$$

\

\noindent \textbf{(S2)} If $d^u < \infty$ then $\theta^p(d^u)>0$ for every $p\in \mathcal{J}$. If $d^u=\infty$ then
$\theta^p(i) > 0$ for every $i\ge 1$ and $p\in \mathcal{J}$.

\



\noindent In the sequel, $vec(A)$ denotes the usual vectorization of an $m\times n$ matrix $A$, i.e., $vec(A) := [a_{11},\ldots,a_{m1},a_{12},\ldots,a_{m2},\ldots, a_{1n},\ldots, a_{mn}]^T$. Here $a_{ij}$ represents the $(i,j)$-th element of a given matrix $A$ and the superscript $T$ denotes the transpose.

\begin{theorem} \label{Conv_LR_infinity}
Assume that a weighted process $U$ satisfies assumptions \textbf{(M1, M2, M3$^{\prime}$, M4$^{\prime}$)} and \textbf{(H1$^{\prime})$}. Then for every sequence of $\mathbb{F}$-stopping times $d^l_{n^\star}$ and $d^u_{n^\star}$ satisfying \textbf{(S1-S2)}, we have

\begin{equation}\label{randomconv}
\sum_{\ell=d^l_{n^\star}}^{d^u_{n^\star}} \xi^{n^\star}(\ell)\rightarrow N\Big(0, \Gamma(d^u)-\Gamma(d^l-1)\Big)~\text{weakly as}~ n^\star\rightarrow \infty.
\end{equation}

\noindent Moreover,

\begin{equation}\label{covconv}
vec\Big(\hat{\Gamma}(n^\star, d^u_{n^\star}) - \hat{\Gamma}(n^\star, d^l_{n^\star}-1)\Big) \rightarrow vec\Big(\Gamma(d^u)-\Gamma(d^l-1)\Big)
\end{equation}

\noindent in probability as $n^\star\rightarrow \infty$.

\end{theorem}

\begin{remark}
Assume $d^u <\infty$ and $\theta^p(d^u) > 0$ for every $p\in \mathcal{J}$ and $\xi^{n^\star}_{m,q}(\ell)$ is square-integrable for every $n^\star\in \mathbb{N}^J, q\in \mathcal{J}$ and $\ell \ge 1$. Then the result stated in Theorem~\ref{Conv_LR_infinity} holds under assumptions~\textbf{(M1-M2-M3-M4)} and \textbf{(H1)}.
\end{remark}

\begin{remark}\label{j<3}
If we have just two populations ($\mathcal{J}=\{1,2\}$) then the underlying covariance structure in Theorem~\ref{Conv_LR_infinity}  simplifies substantially since in this case

$$\langle Q(\ell)a,a\rangle_{\mathbb{R}^J}=\sum_{k=1}^Ja^2_k\phi^2_k(\ell),$$

$$
\langle \hat{Q}(n^\star,\ell) a,a\rangle_{\mathbb{R}^J} =\sum_{k=1}^J a^2_k \hat{\phi}^2_{k,n^\star}(\ell),
$$

\noindent for $a\in \mathbb{R}^J,~n^\star\in \mathbb{N}^J$ and $1\le \ell < \infty$.
\end{remark}

In the sequel, we explore the whole trajectory of the $\mathbb{R}^J$-valued process $\xi^{n^\star}$ weighted by the sequence $\hat{\phi}_{n^\star}:=\{\hat{\phi}_{1,n^\star}(i), \ldots, \hat{\phi}_{J,n^\star}(i)~; i\ge 1 \}$ in a suitable Hilbert space which encodes the quadratic powers of $\xi^{n^\star}$. Doing so, our main motivation and application for the next result is the introduction of the Cram\'er-von Mises test statistics under arbitrary right censoring for infinitely many categories. Let us define the following $\mathbb{R}^J$-valued weighted random field

\begin{equation}\label{glr}
GLR(n^\star, x, r):=\Bigg\{x_q(r)\sum_{\ell=1}^r \sum_{m=1}^{N_J} \xi^{n^{\star}}_{m,q}(\ell);~q=1,\ldots,J\Bigg\}
\end{equation}
for $r\ge 1$, $x=\{x_{q}(i);q=1\ldots, J,~i\ge 1 \}$ and $n^\star\in \mathbb{N}^J$. With this $\mathbb{R}^J$-valued random field, for a given $(n,m)\in \mathbb{N}^2$  with $1\le n\le m$, $x\in \mathbb{R}^{\infty}$ and $n^\star\in \mathbb{N}^J$ we define

\begin{equation}\label{get}
GET(n^\star, x, n ,m):=\big(GLR(n^\star,x ,n),\ldots,GLR(n^\star,x ,m)\big).
\end{equation}

\noindent Of course, under mild assumptions on the weights $(x,U^{n_q}_{n_{q_1}})$ and the sample we can safely embed the process $GET(n^\star, x, \cdot ,\cdot)$ into the Hilbert space $\ell^2(\mathbb{N})$ constituted of square-summable real sequences over $\mathbb{N}$. In the sequel, we make use of the following notation: $\|\cdot\|_{\ell^2}$ stands for the usual norm on the Hilbert space $\ell^2(\mathbb{N})$, $M(\ell):=\text{diag}(\phi_1(\ell),\ldots,\phi_{J}(\ell))$.

In the sequel, for a given $1\le  s \le i < \infty$ we consider the self-adjoint operator $Y(s,i):\mathbb{R}^{Jk(s,i)}\rightarrow \mathbb{R}^{Jk(s,i)}$ defined by the following quadratic form

\begin{eqnarray*}
\langle Y(s , i)a,a\rangle_{\mathbb{R}^{Jk(s,i)}} &:=& \sum_{j=1}^{k(s,i)} \langle M(j+s -1)\Gamma(j+s -1) M(j+s -1)a_j,a_j\rangle_{\mathbb{R}^{J}} \\
& & \\
&+&  \sum_{1\le \ell< j \le k(s,i)} \langle M(\ell + s -1)\Gamma(\ell + s -1) M(j + s-1)a_{\ell},a_j\rangle_{\mathbb{R}^{J}} \\
&+& \sum_{1 \le j < \ell \le k(s,i)} \langle M(\ell + s -1)\Gamma(j + s -1) M(j + s-1)a_{\ell},a_j\rangle_{\mathbb{R}^{J}},
\end{eqnarray*}

\noindent for $a\in \mathbb{R}^{Jk(s,i)}$;~$k(s,i):=i-s+1$. From the definition of the above quadratic form, we notice that for each $s\ge 1$ the restriction of $Y(s, j)$ onto $\mathbb{R}^{Jk(s,i)}$ is equal to $Y(s, i)$ for every $j\ge i$. Therefore, for a given $(s,m)\in \bar{N}^2$ with $1\le s\le m \le \infty$ we shall construct a linear map

\begin{equation}\label{Ycons}
\mathcal{Y}(s,m):\mathbb{R}^\infty\rightarrow \mathbb{R}^\infty
\end{equation}
defined as follows. If $m< \infty$, then we set $\mathcal{Y}(s,m)a:=Y(s,m)(a_1,\ldots, a_{Jk(s,m)})$ so that $\mathcal{Y}(s,m)=Y(s,m)$. If $m=\infty$, then for a given $a\in\mathbb{R}^\infty$ the $Jk(s,i)$-th coordinates of the action $\mathcal{Y}(s,m)a$ is defined by

\begin{equation}\label{actionY}
\{Y(s,i)(a_1,\ldots, a_{Jk(s,i)})\};~ s\le i < \infty.
\end{equation}

\

\begin{theorem} \label{main_asymptotic_result}
Assume that assumptions \textbf{(M1, M2, M3$^{\prime}$, M4$^{\prime}$)} and \textbf{(H1$^{\prime})$} hold and let $1\le d^l< d^u\le \infty$ where $\textbf{(S1-S2)}$ holds. Then the weak limit

\begin{equation}\label{etlimit}
\lim_{n^\star\rightarrow \infty}GET(n^\star, \phi ,d^l,d^u)
\end{equation}

\noindent is a zero-mean Gaussian measure on~$\ell^2$ with covariance operator $\mathcal{Y}(d^l,d^u)$ on $\ell^2$ defined by~(\ref{Ycons}) and~(\ref{actionY}). In particular,

\begin{equation}\label{astat}
\|GET(n^\star, \hat{\phi}_{n^\star} ,d^l_{n^\star},d^u_{n^\star})\|^2_{\ell^2}\rightarrow \sum_{s=1}^{\infty} \sum_{q=1}^{J} \lambda_{sq}\chi^2_{sq}~\text{weakly as}~n^\star\rightarrow \infty,
\end{equation}
where $\{\lambda_{sq} ; s\ge 1, q=1, \ldots , J\}$ are the eigenvalues of $\mathcal{Y}(d^l,d^u)$ and $\{\chi_{sq}^2; s\ge 1, q=1, \ldots , J\}$ is an i.i.d subset of chi-squared random variables with one degree of freedom.
\end{theorem}


The remainder of this paper is devoted to give applications of Theorems~\ref{Conv_LR_infinity} and~\ref{main_asymptotic_result} to the analysis of purely discrete populations typically founded in a lifetime data setting. At first, we exhibit a large class of weighted process which satisfies the assumptions in Theorems~\ref{Conv_LR_infinity} and~\ref{main_asymptotic_result}.

\subsection{Weighted Processes}  \label{weighted}
A large number of weighted processes satisfy \textbf{(M1, M2, M3$^{\prime}$, M4$^{\prime}$)} and \textbf{(H1$^{\prime}$)}. For instance, they can be chosen according to the following generic class

\begin{equation} \label{HFTW}
U^{n_q}_{n_{q_1}} (n^{\star}, \ell ) := \left( \frac{1}{n} \right)^{1/2} u (n^{\star}, \ell) \left(\frac{V^{n_q}(\ell)V^{n_{q_1}}(\ell)}{V^{n^{\star}}(\ell)} \right),~\ell\ge 1,
\end{equation}


\noindent for any pair $q \neq q_1$ in $\mathcal{J}$ and $n^\star=(n_1,\ldots,n_J)\in\mathbb{N}^J$ where we set $n:=\sum_{i=1}^Jn_i$. We assume that the weighted process $u (n^{\star}, \cdot)$ is bounded, it satisfies the measurability assumption \textbf{(M1)} and it converges in probability to a bounded real-valued function $\omega$. A similar class has previously appeared in Andersen et al~\cite{ABGK} for the continuous case. In the sequel, we denote by $\mathcal{K}$ the class of all weighted process which can be represented by (\ref{HFTW}).

\begin{proposition}\label{HAass}
 Let us assume the existence of the limit $b_p = \lim_{n^\star \rightarrow \infty}n_{p} / n$ and $X^p$ is integrable for every $p\in \mathcal{J}$. Then, every weighted process in the class $\mathcal{K}$ satisfies assumptions \textbf{(M1, M2, M3$^{\prime}$, M4$^{\prime}$)} and \textbf{(H1$^{\prime})$}.
\end{proposition}

A significant subclass of $\mathcal{K}$ is given by the classical Tarrone-Ware~\cite{TW} and Harrington-Fleming~\cite{HF} weighted processes. The weighted functionals $u(n^\star,\cdot)$ are given, respectively, by

\begin{equation}\label{hftwfunctions}
 \varphi\left( \frac{V^{n^\star}(\ell)}{n}\right),\quad \Bigg(\frac{\Delta R^{n^\star}(\ell-1)}{V^{n^\star}(\ell-1)} \Bigg)^{\beta} \left( \prod_{j=0}^{\ell-1} \left( 1 - \frac{\Delta R^{n^\star}(j)}{V^{n^\star}(j)}  \right)\right)^\delta
,n^\star\in\mathbb{N}^J,~\ell\ge 1,
\end{equation}

\noindent where $\varphi$ is a bounded continuous function and $\beta$ and $\delta$ are positive constants.

\section{The Log-Rank Statistics}\label{log-ranksection}

In this section, inspired by the weighted log-rank statistics proposed by Gill~\cite{Gill}, Fleming and Harrington~\cite{HF} and Andersen et al~\cite{ABGK}, we propose a test in order to verify the homogeneity of discrete populations in the presence of arbitrary right censoring. Our goal is to derive a class of statistical tests for the null hypothesis

$$H_0 : h^1 (\ell) = h^2(\ell) = \cdots = h^J (\ell); \ell \in\mathbb{N}.$$
Throughout this section, all weighted processes belong to the class $\mathcal{K}$. In the sequel, we denote


\begin{equation}\label{du}
d^u = \sup \left\{ \ell :  \min_{q \in \mathcal{J} } \theta^q (\ell) > 0 \right\},
\end{equation}

\begin{equation}\label{dustar}
d^u_{n^\star}= \sup \left\{ \ell: \min_{q \in \mathcal{J} } V^{n_q}( \ell )  > 0 \right\},~n^\star\in \mathbb{N}^J.
\end{equation}

\noindent One can easily check that $d^u_{n^{\star}}\rightarrow d^u$ in probability as $n^\star\rightarrow \infty$, where $1\le d^u\le\infty$ and $d^u_{n^\star}< \infty$ a.s for every $n^\star\in \mathbb{N}^J$. If $q\in\mathcal{J}$, then we introduce the following general linear $J$-sample statistics

\begin{eqnarray*}
\nonumber LR_{q}({n^{\star}}, j) &:=& \sum_{\ell=1}^{j} \sum_{q_1\neq q}  U^{n_q}_{n_{q_1}} (n^{\star}, \ell ) \left[ \hat{h}^{n_q} (\ell) -  \hat{h}^{n_{q_1}} (\ell)\right] \\
&=& \sum_{\ell=1}^{j} \left( \frac{1}{n} \right)^{1/2} u (n^{\star}, \ell) V^{n_q} (\ell) \left[ \frac{\Delta R^{n_q} (\ell)}{V^{n_q} (\ell)} - \frac{\Delta R^{n^{\star}} (\ell)}{V^{n^{\star}} (\ell)} \right],
\end{eqnarray*}

\noindent for $n^\star\in \mathbb{N}^J,~j\ge 1$. We notice that $LR_{q}(n^\star, \cdot)$ is the $q$-th component of $\xi^{n^\star}$ under the particular null hypothesis $H_0$. Following Theorem~\ref{Conv_LR_infinity} and Proposition~\ref{HAass}, under $H_0$ the  random vector

$$LR({n^{\star}}, d^u_{n^{\star}}) :=(LR_1({n^{\star}},d^u_{n^{\star}}), \ldots , LR_{J}({n^{\star}},d^u_{n^{\star}}))^T$$

\noindent converges weakly to $N\big(0, \Gamma(d^u)\big)$ as $n^\star\rightarrow \infty$, where $\Gamma(d^u)$ admits a consistent estimator $\hat{\Gamma}({n^{\star}},d^u_{n^{\star}})$ given by the matrix induced by the quadratic form~(\ref{covoperatorest}).  Similar to the continuous case, the sum of the components of the random vector $LR$ is null and hence we consider the vector $LR$ without the last component as follows

$$LR_0({n^{\star}},d^u_{n^{\star}}) := (LR_1({n^{\star}},d^u_{n^{\star}}), \ldots , LR_{J-1}({n^{\star}},d^u_{n^{\star}})^T.$$
In what follows, we denote $\Gamma_0(i):= \sum_{\ell=1}^iQ_0(\ell)$ where $Q_0$ is the operator $Q$ defined in~(\ref{covoperator}) without the last row and column. We are now in position to define the weighted log-rank statistics associated to $LR_0$ as follows

\begin{equation}\label{wstatlg}
X^2({n^{\star}}, d^u_{n^{\star}}) := LR_0({n^{\star}}, d^u_{n^{\star}})^T \hat{\Gamma}_0({n^{\star}},d^u_{n^{\star}})^{-1} LR_0({n^{\star}}, d^u_{n^{\star}}),~n^\star \in\mathbb{N}^J,
\end{equation}
where $\hat{\Gamma}_0({n^{\star}},j):=\sum_{\ell=1}^{j}\hat{Q}_0(n^\star, \ell); j\ge 1$. One can easily check that the statistics~(\ref{wstatlg}) is asymptotically chi-square distributed with $J-1$ degrees of freedom, where $\hat{\Gamma}_0({n^{\star}},d^u_{n^{\star}})^{-1}$ is the ordinary inverse.

\begin{remark}\label{obsinv}
Under $H_0$, we denote the intensity function by $h^p(\ell) = h(\ell)$ and $\theta^p (\ell) = \theta(\ell)$ for each $\ell \geq 1$ and  $p \in \mathcal{J}$. It follows from (\ref{covariance_neg}) that the covariance component of $\Gamma$ is negative.  In fact, for each $k \neq r~\text{in}~\mathcal{J}$

\[
\psi(k,r,\ell)= -J ~ \omega^2 (\ell) b_r b_k \theta(\ell) h(\ell)[1-h(\ell)], \quad \ell = d^l, \ldots , d^u.
\]
In this case, we can apply the same arguments of Andersen et al~\cite{ABGK} to conclude that the $rank (\Gamma(\ell))=J-1$ for $d^l\le \ell\le d^u$. By consistency of $\hat{\Gamma}(\ell)$, the probability that $\hat{\Gamma}(\ell)$ has rank $J-1$ increases to unity as $n^{\star}\rightarrow \infty$.

\end{remark}

\begin{remark}\label{obsext}
The class of statistics $X^2({n^{\star}}, \cdot)$ indexed by weighted processes in the class $\mathcal{K}$ is rather general in the sense that it covers the classical cases of Tarone-Ware and Harrington-Fleming but it is not restricted to these cases. In fact, we have proved that one can choose any weighted process $U$ which satisfies the assumptions in Theorem~\ref{Conv_LR_infinity}. Under these assumptions, one can always define a correspondent weighted log-rank statistics of the form~(\ref{wstatlg}).
\end{remark}

\noindent \textit{Type 2 Censoring}. In many practical applications, the censoring scheme is linked to a failure time process. For instance, let us suppose that the stopping time for a product life testing is not fixed a priori, i.e., it is \textit{not} fixed before the beginning of the study but it is chosen later, with the choice influenced by the results of the study up to that time. The so-called Type 2 censoring refers to the case which based on the observed data on that moment, one may want to stop the experiment at some stopping time. Let us suppose the case where the experiment finishes at the occurrence time of the $\beta$-quantile of the observed data $(0 < \beta < 1)$. In order to develop a suitable model for this important type of censoring, we assume $X^p=W^p$ and stop the weighted log-rank process $LR_0 (n^\star , \cdot)$ at the time of the occurrence of the $\beta$-quantile.

Let us consider the following stopping times
$$\ell_{n^{\star}}:= \inf \{\ell\ge 1 ; R^{n^{\star}} (\ell) / n \geq \beta \},$$

$$\ell^0:= \inf \{\ell : b_1 F^1(\ell) + \ldots + b_J F^J (\ell) \geq \beta\},$$
where $F^q$ is the distribution function of $W^q,~q\in \mathcal{J}$. One can easily check that the following convergence holds
$$\ell _{n^{\star}}\rightarrow \ell^0,$$
in probability as $n^\star\rightarrow \infty$. As a consequence of Theorem~\ref{Conv_LR_infinity} and Proposition~\ref{HAass}, we arrive at the following result.

\begin{corollary} \label{Type2_Conv}
Under $H_0$ the log-rank statistics  $LR_0({n^{\star}},\ell_{n^{\star}})$ converges weakly to $N\big(0, \Gamma_0(\ell^0)\big)$ as $n^\star\rightarrow \infty$, where $\Gamma_0(\ell^0)$ admits a consistent estimator $\hat{\Gamma}_0({n^{\star}}, \ell_{n^{\star}})$.
 \end{corollary}

\section{The Cram\'er-von Mises Statistics}\label{CVMsection}

The goal of this section is to propose a class of statistical tests for the null hypothesis $H_0$ which is consistent to any alternative hypothesis. As a consequence of Theorem~\ref{main_asymptotic_result}, we are able to introduce a Cram\'er-von Mises statistics for purely discrete populations under censoring as described in Section~\ref{multmodel}. In the remainder of this section, we encode the first observed categories by the following sequence of $\mathbb{F}$-stopping times

\begin{equation}\label{dlstar}
d^l_{n^{\star}} := \inf \{ \ell : \Delta R^{n^{\star}} (\ell ) > 0 \}
\end{equation}
and we assume that all weighted processes belong to the class $\mathcal{K}$. One can easily check that

$$
d^l_{n^{\star}}\rightarrow  d^l
$$
in probability as $n^\star\rightarrow \infty$ where

\begin{equation}\label{dl}
d^l:= \inf \{ \ell : b_1 h^1(\ell) + \cdots +b_J h^J (\ell) > 0\}
\end{equation}
for $b_p = \lim_{n^\star \rightarrow \infty}n_{p} / n; p\in \mathcal{J}$. Let us now introduce a version of the Cram\'er-von Mises statistics in order to test homogeneity of discrete populations in the presence of arbitrary right censoring with infinitely many categories. From Theorem~\ref{main_asymptotic_result}, we only need to consider the particular case when the $q$-th component of $\xi^{n^\star}$ is the log-rank statistics with dimension $J-1$. That is, $GET(n^\star, \hat{\phi}_{n^\star},i, j)$ in~(\ref{get}) is composed by the following weighted random field

\begin{equation}\label{cvmcomp}
GLR(n^\star, \hat{\phi}_{n^\star} ,r)=\{\hat{\phi}_{q, n^\star}(r)LR_{q}(n^\star,r);q=1,\ldots,J-1\}
\end{equation}
for $r\ge 1, ~n^\star\in\mathbb{N}^{J}.$ We set $M_0(\ell):=\text{diag}(\phi_1(\ell),\ldots,\phi_{J-1}(\ell))$ and the covariance operator in Theorem~\ref{main_asymptotic_result} (see~(\ref{Ycons})) is defined in the same way but $M_0$ and $\Gamma_0$ instead of $M$ and $\Gamma$ in the quadratic form. We denote this linear operator by $\mathcal{Y}_0(d^l,d^u)$ so that the $(J-1)k(d^l,i)$-th coordinates of the action $\mathcal{Y}_0(d^l,d^u)a$ is $\{Y_0(d^l,i)(a_1,\ldots,a_{(J-1)k(d^l,i)})\}$ for $a\in \ell^2$ and $k(d^l,i)=i-d^l+1; i\ge d^l$.

The natural candidate for the estimator of the operator $\mathcal{Y}_0(d^l,d^u):\ell^2\rightarrow \ell^2$ can be constructed in a natural way as follows. To shorten notation, we introduce the random set $L(d^l_{n^\star},d^u_{n^\star}) = \{ d^l_{n^\star} \leq \ell \leq d^u_{n^\star}: \Delta R^{n^\star} (\ell) > 0  \}$ of observable categories and we denote $L(n^{\star})$ its cardinality.   For a given $n^\star\in\mathbb{N}^J$ and $a\in \ell^2$, we define the action

$$\hat{\mathcal{{Y}}}_0(d^l_{n^\star}, d^u_{n^\star})a$$ as the real sequence where the $(J-1)L(n^{\star})$-th coordinates are given by
$$\hat{Y}_0(d^l_{n^\star}, d^u_{n^\star})(a_1\ldots, a_{(J-1)L(n^{\star})})$$
and $\hat{Y}_0(d^l_{n^\star},d^u_{n^\star})$ is the self-adjoint random operator defined by the following quadratic form over $\mathbb{R}^{(J-1)L(n^{\star})}$

\begin{eqnarray*}
\langle \hat{Y}_0(d^l_{n^{\star}},d^u_{n^\star})a,a\rangle &=& \sum_{j \in L(d^l_{n^\star},d^u_{n^\star})} \langle \hat{M}_0(j)\hat{\Gamma}_0(j) \hat{M}_0(j)a_j,a_j\rangle_{\mathbb{R}^{J-1}}\\
& &\\
&+&\sum_{\{ \ell< j : \ell, j \in L(d^l_{n^\star},d^u_{n^\star})\}} \langle \hat{M}_0(\ell) \hat{\Gamma}_0(\ell) \hat{M}_0(j)a_{\ell},a_j\rangle_{\mathbb{R}^{J-1}}\\
& &\\
&+&\sum_{ \{ j < \ell : \ell, j \in L(d^l_{n^\star},d^u_{n^\star}) \}} \langle \hat{M}_0(\ell) \hat{\Gamma}_0(j) \hat{M}_0(j)a_{\ell},a_j\rangle_{\mathbb{R}^{J-1}},
\end{eqnarray*}
where $\hat{M}_0(\cdot):=\text{diag}(\hat{\phi}_{1,n^\star}(\cdot),\ldots,\hat{\phi}_{J-1,n^\star}(\cdot))$ and $a\in \mathbb{R}^{(J-1)L(n^{\star})}$. Therefore, $\mathcal{\hat{Y}}_0(d^l_{n^\star}, d^u_{n^\star}):\ell^2\rightarrow \ell^2$ is a well-defined sequence of self-adjoint finite-rank random operators.

In view of the log-rank composition in~(\ref{cvmcomp}), we are now in position to introduce the Cram\'er-von Mises statistics associated to the general discrete censoring model described in Section~\ref{multmodel} as follows

$$
 CVM (n^\star,d^l_{n^{\star}} , d^u_{n^{\star}}) := \|GET(n^\star, \hat{\phi}_{n^\star}, d^l_{n^{\star}},d^u_{n^{\star}})\|^2_{\ell^2}; n^\star\in \mathbb{N}^{J-1}.
$$
As a consequence of Theorem~\ref{main_asymptotic_result} and Proposition~\ref{HAass}, we arrive at the following result.
\begin{theorem} \label{main_CVM}
Assume that $U$ belongs to the class $\mathcal{K}$, the growth condition in Proposition~\ref{HAass} holds and let
$(d^l,d^u,d^l_{n^\star},d^u_{n^\star})$ be the categories and the stopping times defined by~(\ref{dl}),~(\ref{du}),~(\ref{dlstar}) and~(\ref{dustar}), respectively. Then, under $H_0$

\begin{equation}\label{BCVM}
CVM (n^\star,d^l_{n^{\star}} , d^u_{n^{\star}})\rightarrow \sum_{s=1}^{\infty} \sum_{q=1}^{J-1} \lambda_{sq}\chi^2_{sq}~\text{weakly as}~n^\star\rightarrow \infty,
\end{equation}
\noindent where $\{\lambda_{sq} ; s\ge 1, q=1, \ldots , J-1\}$ are the eigenvalues of the covariance operator $\mathcal{Y}_0(d^l,d^u)$. In particular, if $X^q$ is square-integrable for every $q\in \mathcal{J}$ then

\begin{equation} \label{apr_CVM}
 \Lambda(n^{\star}):=\sum_{s=1}^{L(n^{\star})} \sum_{q=1}^{J-1} \hat{\lambda}_{sq} \chi_{sq}^2  1\!\!1_{ A(n^{\star})}\rightarrow  \sum_{s=1}^{\infty} \sum_{q=1}^{J-1} \lambda_{sq}\chi^2_{sq}
\end{equation}

\noindent weakly as~$n^\star\rightarrow \infty$, where $\{\hat{\lambda}_{sq} : 1\le s\le L(n^{\star}), q=1, \ldots , J-1\}$ are the random eigenvalues of the covariance operator estimator $\hat{\mathcal{Y}}_0( d^l_{n^{\star}},d^u_{n^{\star}})$ and

\[
A(n^{\star}) := \left\{ \hat{\mathcal{Y}}_0( d^l_{n^\star}, d_{n^{\star}}^u)~\text{is non-negative}\right\},
\]
so that $\mathbb{P}\big(A(n^\star)\big)\rightarrow 1$ as $n^\star\rightarrow \infty$.
\end{theorem}
We notice from~(\ref{apr_CVM}) that the $P$-value for the hypothesis test $H_0$ is given by $\mathbb{P} [ \Lambda(n^{\star}) > CVM (n^\star,d^l_{n^{\star}} , d^u_{n^{\star}}) \mid H_0]$. The approximate law $\Lambda(n^{\star})$ is a weighted sum of independent chi-squared random variables and hence several algorithms to evaluate the $P$-value are available. See e.g Duchesne and Micheaux~\cite{DUMI} for a recent discussion.

%

\section{Simulation}\label{simsection}
In this section, we perform a simple simulation study to evaluate the behavior of the classes of statistics proposed in this paper.  We analyze the effect of the sample size, the proportion of censored data and the number of populations. Here, we assume the variable of interest follows a Poisson distribution with parameter $\lambda_p$ for any $p\in\mathcal{J}$ and the censoring variable also follow a Poisson distribution with parameter $\lambda_c$.  In this section, all statistics are considered in terms of of the weight $u(n^{\star}, \ell) =1$ in the class $\mathcal{K}$. The goal is to test the hypothesis $H_0 : \lambda_1 = \lambda_2 = \cdots = \lambda_J$ which is equivalent to the hypothesis $H_0 : h^1(\ell) = h^2(\ell) = \cdots = h^J(\ell);~\ell \geq 1$. The simulation is performed by means of the software R (\cite{Software_R}) and the $P$-value is evaluated via Davies algorithm (\cite{Davies}).

In order to evaluate the convergence of the proposed statistics, we have sampled from several populations with Poisson distribution with $\lambda_1 =100$ and taking into account different sample sizes and censoring variables. It was generated $10000$ samples for each sample size (SS), $\lambda_c$ (without censoring, 90 and 100) and populations (2,4 and 8). Tables \ref{WCV}, \ref{100PoissonCensoring} and \ref{90PoissonCensoring} show the empirical significance level related to the nominal significance level $\alpha=0.05$.


\begin{table}[h]
  \centering
\begin{tabular}{|r|c|c|c|c|c|c|}
 \hline

           & \multicolumn{ 2}{c|}{Two Populations} & \multicolumn{ 2}{c|}{Four Populations} & \multicolumn{ 2}{c|}{Eight Populations} \\ \hline

      SS &        CVM &         LR &        CVM &         LR &        CVM &         LR \\ \hline

        50 &     0.0493 &     0.0564 &     0.0478 &     0.0569 &     0.0706 &      0.072 \\ \hline

       100 &     0.0492 &     0.0531 &     0.0495 &     0.0516 &     0.0652 &     0.0602 \\ \hline

       150 &      0.054 &     0.0549 &     0.0526 &     0.0501 &     0.0599 &     0.0593 \\ \hline

       200 &     0.0478 &     0.0493 &     0.0536 &     0.0527 &     0.0595 &     0.0544 \\ \hline

       250 &      0.051 &     0.0524 &     0.0526 &     0.0528 &     0.0549 &     0.0563 \\ \hline

       300 &     0.0493 &     0.0528 &     0.0534 &     0.0523 &     0.0545 &     0.0541 \\ \hline

\end{tabular}

\caption{Without Censoring.}\label{WCV}
\end{table}

\begin{table}[h]
  \centering
   \begin{tabular}{|r|c|c|c|c|c|c|}
   \hline

           & \multicolumn{ 2}{c|}{Two Populations} & \multicolumn{ 2}{c|}{Four Populations} & \multicolumn{ 2}{c|}{Eight Populations} \\ \hline

        SS &        CVM &         LR &        CVM &         LR &        CVM &         LR \\ \hline

       50 &     0.0506 &     0.0513 &     0.0457 &     0.0588 &     0.0535 &     0.0682 \\ \hline

       100 &     0.0534 &     0.0507 &     0.0454 &     0.0538 &     0.0512 &     0.0595 \\ \hline

       150 &     0.0526 &     0.0504 &     0.0497 &     0.0558 &     0.0535 &     0.0581 \\ \hline

       200 &     0.0504 &     0.0519 &     0.0482 &     0.0525 &     0.0529 &     0.0558 \\ \hline

       250 &     0.0475 &     0.0487 &     0.0504 &     0.0525 &     0.0456 &      0.051 \\ \hline

       300 &     0.0517 &     0.0547 &     0.0482 &     0.0499 &     0.0505 &     0.0539 \\ \hline

\end{tabular}
 \caption{Censoring Variable: Poisson with $\lambda_c =100$.}\label{100PoissonCensoring}
\end{table}

\begin{table}[h]
  \centering

  \begin{tabular}{|r|c|c|c|c|c|c|}
  \hline

           & \multicolumn{ 2}{c|}{Two Populations} & \multicolumn{ 2}{c|}{Four Populations} & \multicolumn{ 2}{c|}{Eight Populations} \\ \hline

        SS &        CVM &         LR &        CVM &         LR &        CVM &         LR \\ \hline

        50 &     0.0491 &     0.0477 &     0.0285 &     0.0711 &     0.0326 &     0.1051 \\ \hline

       100 &     0.0483 &     0.0484 &      0.038 &     0.0578 &     0.0387 &      0.072 \\ \hline

       150 &      0.048 &     0.0485 &     0.0457 &     0.0531 &     0.0446 &     0.0575 \\ \hline

       200 &     0.0474 &     0.0482 &      0.043 &     0.0552 &     0.0451 &     0.0635 \\ \hline

       250 &     0.0516 &     0.0532 &     0.0441 &     0.0506 &     0.0425 &     0.0569 \\ \hline

       300 &     0.0481 &      0.048 &     0.0429 &     0.0532 &     0.0451 &     0.0591 \\ \hline

\end{tabular}
  \caption{Censoring Variable: Poisson with $\lambda_c =90$.}\label{90PoissonCensoring}
\end{table}

In all cases, the simulation shows that both test statistics have approximately the same behavior even for small and moderate sample sizes.  Unsurprisingly, the number of populations involved in the analysis affects the convergence of the log-rank and Cram\'er-von Mises statistics. Moreover, as the proportion of censored data increases, the number of populations becomes more relevant for the convergence.

As pointed out in the Introduction, the logrank test has little power for crossing hazard functions. In order to evaluate the performance of the Cram\'er-von Mises test developed in this paper, we analyze one classical discretely recorded data set with crossing intensities described in Klein and Moeschberger (\cite{KM}, pp. 211). A clinical trial of chemotherapy against chemotherapy combined with radiotherapy in the treatment of locally unresectable gastric cancer was conducted by the Gastrointestinal Tumor Study Group. In this trial, forty-five patients were randomly divided into two groups and medically accompanied for eight years. We wish to test the null hypothesis $H_0$ that the intensity functions of the two groups are the same by using the Cram\'er-von Mises statistics developed in this paper. By setting $(u(n^{\star}, \cdot) = 1)$ in~(\ref{HFTW}), we obtain $CVM=0.0926$ with $P$-value of $0.029$, so the null hypothesis of no difference among intensity functions between the groups is rejected at the $5\%$ level. The same hypothesis test based on the continuous versions of Renyi ($P$-value = 0.053) and Cram\'er-von Mises ($P$-value = 0.06) statistics do not reject $H_0$ at the same $5\%$ level. This result stresses the importance of modeling discrete data with methodologies based on purely discrete distributions. In this particular case, classical methods based on continuous distributions fail to reject $H_0$ while the discrete Cram\'er-von Mises test developed in this article successfully reject it at a given level of significance.

\section{Proof of Theorem~\ref{Conv_LR_infinity}}\label{proofTH1}

In this section, we provide the proof of the asymptotic result stated in Theorem~\ref{Conv_LR_infinity}. Proofs of Lemmas~\ref{multimardiff},~\ref{con_ell_2},~\ref{mainres1} and Proposition~\ref{multlema} are given in the Appendix~\ref{appendixlogrank}. The following simple remark gives the expected probabilistic structure of $\xi^{n^*}_{m,q}$ given by~(\ref{marraydif}).

\begin{lemma}\label{multimardiff}
Assume that $U$ satisfies assumption~\textbf{(M1)} and $\{\xi^{n^\star}_{m,q}(\ell); 1\le m\le N_J, \ell\ge 1, n^\star\in \mathbb{N}^J, q\in\mathcal{J} \}$ is a subset of $L^2(\mathbb{P})$. Then for each $\ell\ge 1$, $n^\star\in \mathbb{N}^J$ and $q\in \mathcal{J}$, $\{ (\xi^{n^{\star}}_{m,q}(\ell); 1\le m \le N_J\}$ is a square-integrable martingale-difference w.r.t the filtration $\mathcal{G}^{n^\star} (\ell) = \{\mathcal{G}^{n^\star}_{m}( \ell); 1\le m \le N_J\}$.
\end{lemma}

Throughout this section, we assume that $\mathbb{E}|\xi^{n^\star}_{m,q}(\ell)|^2< \infty$ for every $ n^\star\in \mathbb{N}^J, 1\le m\le N_J, \ell \ge 1,$ and $q\in\mathcal{J}$.

\begin{lemma} \label{con_ell_2}
Assume that a weighted process satisfies assumptions \textbf{(M1, M2, M3, M4)}. Let $i$ be a positive integer such that $\theta^p(i)>0$ for every $p\in \mathcal{J}$. Then for each $\ell\in \{1,\ldots,i \}$ and $q\in \mathcal{J}$,
\begin{equation}\label{w1}
 \sum_{m=1}^{N_J} \xi^{n^{\star}}_{m,q}(\ell)\rightarrow  N\Big(0, \phi^2_q(\ell)\Big)\quad \text{weakly as}~n^\star\rightarrow \infty.
\end{equation}

\noindent The asymptotic variance $\phi^2_q(\ell)$ is the limit in probability of $\hat{\phi}^2_{q,n^\star}(\ell)$ as $n^\star\rightarrow \infty$ for each $\ell\in \{1,\ldots, i\}$.
\end{lemma}

\begin{lemma}\label{mainres1}
Assume that a weighted process $U$ satisfies assumptions \textbf{(M1, M2, M3, M4)}. Let $i$ be a positive integer such that $\theta^p(i)>0$ for every $p\in \mathcal{J}$. Then for each $q\in \mathcal{J}$, the random variables $\{\sum_{m=1}^{N_J}\xi^{n^\star}_{m,q}(\ell);1\le \ell \le i \}$ are asymptotically independent and

$$\sum_{\ell=1}^i  \sum_{m=1}^{N_J} \xi^{n^{\star}}_{m,q}(\ell)\rightarrow N\Big(0,\sum_{\ell=1}^i\phi^2_q(\ell)\Big)\quad \text{weakly as}~n^\star\rightarrow \infty.$$
\end{lemma}

\begin{proposition}\label{multlema}
Assume that a weighted process $U$ satisfies assumptions \textbf{(M1-M2-M3-M4)} and \textbf{(H1)}. Let $i$ be a positive integer such that $\theta^p(i)>0$ for every $p\in \mathcal{J}$. Then for each $\ell \in \{1,\ldots,i \}$

$$\xi^{n^\star}(\ell)\rightarrow N\big(0,Q(\ell)\big)$$


\noindent weakly as $n^\star\rightarrow \infty$, where

\begin{equation}\label{qconv}
vec\Big(\hat{Q}(n^\star,\ell)\Big)\rightarrow vec\Big(Q(\ell)\Big)
\end{equation}

\noindent in probability as $n^\star\rightarrow \infty$ for each $\ell\in\{1,\ldots,i \}$. In particular,

\begin{equation}\label{indvec}
\sum_{\ell=1}^i\xi^{n^\star}(\ell)\rightarrow N(0,\Gamma(i))~\quad \text{weakly as}~n^\star\rightarrow \infty.
\end{equation}

\end{proposition}

\

\noindent \textit{Proof of Theorem~\ref{Conv_LR_infinity}.}

\noindent Throughout this proof $C$ is a constant which may defer from line to line. At first, we assume that $d^l=1, d^u=\infty$ and $\theta^p(i)>0$ for every $p\in \mathcal{J}$ and positive integer $i\ge 1$. By observing the identity~(\ref{imp1}) together with~\textbf{(M3$^{\prime}$-M4$^{\prime}$)} we readily see that $\xi^{n^\star}_{m,q}(\ell)$ is square-integrable for every $q\in\mathcal{J},\ell\ge 1, n^\star\in\mathbb{N}^J$ and $m=1\ldots, N_J$ so we are able to apply Proposition~\ref{multlema} accordingly. Let $\{\tilde{W}(i)=(W_1(i),\ldots, W_J(i)); i\ge 1\}$ be the $\mathbb{R}^J$-valued $\mathbb{F}$-martingale with independent increments given in the proof of Proposition~\ref{multlema}. We claim that

\

(a) $\Gamma(\infty)=\sum_{\ell=1}^\infty Q(\ell)$ is a well-defined nonnegative self-adjoint operator on $\mathbb{R}^J$.

(b) The weak convergence holds $\sum_{\ell=1}^\infty\xi^{n^\star}(\ell)\rightarrow N(0,\Gamma(\infty))$ as $n^\star\rightarrow \infty$.

\

\noindent Let us check (a). For a given $a\in \mathbb{R}^J$, we know that $\langle \tilde{W}(i), a\rangle_{\mathbb{R}^J}$ has Gaussian law $N\big(0,\sum_{\ell=1}^i\langle Q(\ell)a, a\rangle\big)$ for each $i\ge 1$. The definition of the covariance operator $Q$ in~(\ref{covoperator}) and assumptions
\textbf{(M3$^{\prime}$-M4$^{\prime}$)} and \textbf{(H1$^{\prime}$)} yield $\sup_{i\ge 1}\sum_{\ell=1}^i\|Q(\ell)\| < \infty$. In particular, the following estimate holds

\[
\sum_{\ell=1}^i\langle Q(\ell)a,a\rangle_{\mathbb{R}^J} \leq \sum_{\ell=1}^{\infty} \sum_{k=1}^J a^2_k \left(\sum_{q_1\neq k}\alpha_{k,q_1}^{q_1}(\ell) + \sum_{q_1\neq k}\alpha_{k,q_1}^{k}(\ell)+ 2\sum_{(q_1,q_2)\in A_k} |\beta^k_{q_1,q_2}(\ell)| \right) +
\]

\[
2\sum_{\ell=1}^{\infty} \sum_{1\le r< k \le J} |a_ra_k| \left(\sum_{q_1 \in A(k,r)}|\gamma^{k,r}_{q_1}(\ell)|+\sum_{q_1\neq k}|\eta^{k,r}_{q_1}(\ell)| + \sum_{q_2\neq r}|\eta^{r,k}_{q_2}(\ell)| \right) < \infty,
\]
\noindent for every $i \ge 1$ and $a\in \mathbb{R}^J$. Hence, $\Gamma(\infty):=\sum_{\ell=1}^\infty Q(\ell)$ converges absolutely on the space of matrices and it is the self-adjoint non-negative operator associated to the quadratic form $\sum_{\ell=1}^\infty\langle Q(\ell)a,a\rangle_{\mathbb{R}^J};~a\in \mathbb{R}^J.$ Now let us check (b). From Proposition~\ref{multlema} we know that $\tilde{W}(i)\sim N(0,\Gamma(i))$ and hence the previous argument allows us to define the $N\big(0,\Gamma(\infty)\big)$-Gaussian variable $\tilde{W}(\infty):=\lim_{i\rightarrow \infty}\tilde{W}(i)$ (weak sense).

To shorten  notation, we set $\tilde{W}^{n^\star}(i):=\sum_{\ell=1}^i\xi^{n^\star}(\ell); i\ge 1,n^\star\in \mathbb{N}^J$. By using the same argument as in~(\ref{zmart}) one can easily check that $\{\sum_{\ell=1}^i\xi^{n^\star}_{m,q}(\ell); i\ge 1\}$ is an  $\mathbb{F}$-martingale-difference so that $\{\langle \tilde{W}^{n^\star}(i), a\rangle_{\mathbb{R}^J}; i\ge 1 \}$ is an $\mathbb{F}$-martingale for every $a\in \mathbb{R}^J$. For each $i\ge 1$, Proposition~\ref{multlema} yields $\tilde{W}^{n^\star}(i)\rightarrow \tilde{W}(i)\sim N(0,\Gamma(i))$ weakly as $n^\star\rightarrow \infty$. We now claim that the weak limit $\tilde{W}^{n^\star}(\infty):=\lim_{i\rightarrow \infty}\tilde{W}^{n^\star}(i)$ exists for each $n^\star\in \mathbb{N}^J$. In fact, for a given $a\in \mathbb{R}^J$ and $n^\star\in \mathbb{N}^J$, a straightforward but lengthy calculation shows that the quadratic variation of the martingale $\langle \tilde{W}^{n^\star}(\cdot), a\rangle_{\mathbb{R}^J}$ at a given point $i\ge 1$ can be estimated as follows

\[
\sum_{\ell=1}^i|\Delta\langle \tilde{W}^{n^\star}(\ell),a\rangle|^2\leq \sum_{\ell=1}^{\infty} \sum_{k=1}^J a^2_k \left( \sum_{q_1\neq q} \frac{|U^{n_q}_{n_{q_1}}(n^\star,\ell)|^2}{V^{n_{q_1}}(\ell)}+ \right.
\]

\[
\left.  \sum_{q_1\neq q} \frac{|U^{n_q}_{n_{q_1}}(n^\star,\ell)|^2}{V^{n_{q}}(\ell)} +
 2\sum_{(q_1,q_2)\in A_q}\frac{|U^{n_q}_{n_{q_1}}(n^\star,\ell)U^{n_q}_{n_{q_2}}(n^\star,\ell)|}{V^{n_q}(\ell)} \right) +
\]

\[
\sum_{\ell=1}^{\infty} 2\sum_{1\le r<k\le J} |a_ra_k| \left(\sum_{q_1\neq k}\frac{|U^{n_k}_{n_{q_1}}(n^\star,\ell) U^{n_r}_{n_{q_1}}(n^\star,\ell)|}{V^{n_{q_1}}(\ell)} + \right.
\]

\[ \left. \sum_{q_1\neq k}\frac{|U^{n_k}_{n_{q_1}}(n^\star,\ell) U^{n_r}_{n_{k}}(n^\star,\ell)|}{V^{n_{k}}(\ell)} +
\sum_{q_2\neq r}\frac{|U^{n_k}_{n_{r}}(n^\star,\ell) U^{n_r}_{n_{q_2}}(n^\star,\ell)|}{V^{n_{r}}(\ell)} \right).
\]
Therefore, we may use assumptions~\textbf{(M3$^{\prime}$-M4$^{\prime}$-H1$^{\prime}$)} to ensure that for each $a\in \mathbb{R}^J$

\begin{equation}\label{mma1}
\sup_{i\ge 1,n^\star\in \mathbb{N}^J}\mathbb{E}\sum_{\ell=1}^i|\Delta\langle \tilde{W}^{n^\star}(\ell),a\rangle|^2 < \infty
\end{equation}
and hence the Doob maximal inequality jointly with the Martingale convergence theorem yield $\lim_{i\rightarrow \infty}\langle \tilde{W}^{n^\star}(i),a \rangle_{\mathbb{R}^J} $ in probability for every $a\in\mathbb{R}^J$ and $n^\star\in \mathbb{N}^J$. This fact together with Cramer-Wold allow us to conclude that $\tilde{W}^{n^\star}(i)\rightarrow \sum_{\ell=1}^\infty \xi^{n^\star}(\ell)$ weakly on $\mathbb{R}^J$ as $i\rightarrow \infty$ for each $n^\star\in \mathbb{N}^J$. At this point, we know that
\begin{equation}\label{iterated}
\lim_{i\rightarrow \infty}\lim_{n^\star\rightarrow \infty}\tilde{W}^{n^\star}(i)=\tilde{W}(\infty).
\end{equation}
Fix $a\in\mathbb{R}^J$ and denote $\tilde{W}^{n^\star}(\infty)=\sum_{\ell=1}^\infty \xi^{n^\star}(\ell); n^\star\in \mathbb{N}^J$. Since $\{\langle \tilde{W}^{n^\star}(i),a \rangle_{\mathbb{R}^J}; 1\le i\le \infty \}$ is a closed martingale then we can estimate the $L^2$-norm from the quadratic variation as follows. For a given $\varepsilon>0$, the following estimate holds for each $n^\star\in \mathbb{N}^J$ and~$i\ge 1$ as follows

\begin{eqnarray}
\nonumber\mathbb{P}\Big\{\big|\langle \tilde{W}^{n^\star}(\infty)-\tilde{W}^{n^\star}(i),a\rangle_{\mathbb{R}^J}  \nonumber \big|> \varepsilon\Big\}&\le& \frac{1}{\varepsilon^2}\mathbb{E}\Big|\big \langle \langle \nonumber\tilde{W}^{n^\star}(\cdot),a\rangle_{\mathbb{R}^J}  \big \rangle(\infty) - \big\langle \langle \nonumber\tilde{W}^{n^\star}(\cdot),a\rangle_{\mathbb{R}^J}  \big \rangle(i)\Big|\\
\nonumber& &\\
\nonumber&=&\frac{1}{\varepsilon^2}\mathbb{E}\sum_{\ell=i+1}^\infty\Delta \big \langle \langle  \nonumber\tilde{W}^{n^\star}(\cdot),a   \rangle_{\mathbb{R}^J}  \big \rangle(\ell)\\
\nonumber& &\\
\label{deltaarg}&=& \frac{1}{\varepsilon^2}\mathbb{E}\sum_{\ell=i+1}^\infty|\langle  \tilde{W}^{n^\star}(\ell)-\tilde{W}^{n^\star}(\ell-1),a   \rangle_{\mathbb{R}^J}|^2.
\end{eqnarray}
Assumptions~\textbf{(M3$^{\prime}$-M4$^{\prime}$-H1$^{\prime}$)} (in particular~(\ref{mma1})) then imply that

\begin{equation}\label{state1}
\lim_{i\rightarrow \infty}\limsup_{n^\star\rightarrow \infty}\mathbb{P}\Big\{\big|\langle \tilde{W}^{n^\star}(\infty)-\tilde{W}^{n^\star}(i),a\rangle_{\mathbb{R}^J}  \big|> \varepsilon\Big\}=0.
\end{equation}
From [\cite{Billingsley}, Theorem 3.2~page 28] and Cramer-Wold, the convergence~(\ref{state1}) allows us to exchange the limits in~(\ref{iterated}) and therefore $\tilde{W}^{n^{\star}} (\infty) \rightarrow \tilde{W} (\infty)$ weakly as $n^\star\rightarrow \infty$, thus we conclude (b). In particular, we have shown

\begin{equation}\label{pre22}
\tilde{W}^{n^\star}(d^u)-\tilde{W}^{n^\star}(d^l-1)\rightarrow N(0,\Gamma(d^u)-\Gamma(d^l-1))
\end{equation}
weakly as $n^\star\rightarrow \infty$. We now proceed by using the above arguments on the set $\{d^l , \ldots , d^u\}$ and this time we have to play with the stopping times $d^l_{n^\star}$ and $d^u_{n^\star}$. We fix $a\in \mathbb{R}^J$ and we write

\[
\langle \tilde{W}^{n^\star}(d^u_{n^{\star}})-\tilde{W}^{n^\star}(d^l_{n^{\star}}),a\rangle_{\mathbb{R}^J} = \langle \tilde{W}^{n^\star}(d^u)-\tilde{W}^{n^\star}(d^l),a\rangle_{\mathbb{R}^J} +
\]

\[
\langle \tilde{W}^{n^\star}(d^u_{n^{\star}})-\tilde{W}^{n^\star}(d^u),a\rangle_{\mathbb{R}^J} + \langle \tilde{W}^{n^\star}(d^l)-\tilde{W}^{n^\star}(d^l_{n^{\star}}),a\rangle_{\mathbb{R}^J}.
\]
By considering $ c^u_1 (n^{\star}) = d^u \wedge d^u_{n^{\star}}$ and $ c^u_2 (n^{\star}) = d^u \vee d^u_{n^{\star}}$, we may follow the same steps given in~(\ref{mma1}) and (\ref{deltaarg}) to show that for every $\varepsilon>0$

$$\mathbb{P}[|\langle \tilde{W}^{n^\star}(d^u)-\tilde{W}^{n^\star}(d^u_{n^\star}), a\rangle_{\mathbb{R}^J}|> \varepsilon]\le \frac{1}{\varepsilon^2}\mathbb{E}\sum_{\ell=c^u_1(n^\star)+1}^{c^u_2(n^\star)}|\Delta\langle \tilde{W}^{n^\star}(\ell),a \rangle_{\mathbb{R}^J}|^2\rightarrow 0$$
as $n^\star\rightarrow \infty$. The same argument also applies to $d^l \wedge d^l_{n^{\star}}$ and $d^l \vee d^l_{n^{\star}}$ and therefore Cramer-Wold and~(\ref{pre22}) allow us to conclude that (\ref{randomconv}) holds.

The limit~(\ref{covconv}) when $d^u< \infty$ is a direct consequence of Proposition~\ref{multlema} and in particular~(\ref{qconv}). So we only need to prove the case $d^u=\infty$, i.e., $vec\Big(\hat{\Gamma}(n^\star, d^u_{n^\star})\Big)\rightarrow vec\Big(\Gamma(\infty)\Big)$ in probability as $n^\star\rightarrow \infty$. In fact, one has to check the following convergence in probability

\begin{equation}\label{ccc1}
\sum_{\ell=1}^{d^u_{n^\star}}\hat{\phi}^2_{q,n^\star}(\ell)\rightarrow \sum_{\ell=1}^\infty\phi^2_q(\ell); q\in \mathcal{J};
\end{equation}

\begin{equation}\label{ccc2}
\sum_{\ell=1}^{d^u_{n^\star}}\hat{\psi}_{n^\star}(k,r,\ell)\rightarrow \sum_{\ell=1}^\infty\psi(k,r,\ell);~1\le r < k \le J.
\end{equation}
as $n^\star\rightarrow \infty$. In the sequel, to shorten notation we write $\eta^p(\ell)=h^p(\ell)[1-h^p(\ell)]$ for $p\in \mathcal{J},\ell\ge 1$ and we proceed componentwise. For a given $q_1\neq q$ in $\mathcal{J}$, we shall write

\begin{eqnarray}
\nonumber\sum_{\ell=1}^{d^u_{n^\star}}\frac{|U^{n_{q}}_{n_{q_1}}(n^\star,\ell)|^2}{V^{n_{q_1}}(\ell)}\hat{\eta}^{n_{q_1}}(\ell)
\nonumber-\sum_{\ell=1}^\infty\alpha^{q_1}_{q,q_1}(\ell)\eta^{q_1}(\ell)&=& \nonumber\sum_{\ell=1}^{d^u_{n^\star}}\frac{|U^{n_{q}}_{n_{q_1}}(n^\star,\ell)|^2}{V^{n_{q_1}}(\ell)}
\nonumber[\hat{\eta}^{n_{q_1}}(\ell)-\eta^{q_1}(\ell)]\\
\nonumber& &\\
\nonumber&+&\sum_{\ell=1}^{d^u_{n^\star}}\Bigg\{\frac{|U^{n_{q}}_{n_{q_1}}(n^\star,\ell)|^2}{V^{n_{q_1}}(\ell)}
\nonumber-\alpha^{q_1}_{q,q_1}(\ell)\Bigg\}\eta^{q_1}(\ell)\\
\nonumber& &\\
\nonumber&-& \sum_{\ell=d^u_{n^\star}+1}^\infty\alpha^{q_1}_{q,q_1}(\ell)\eta^{q_1}(\ell)\\
\nonumber& &\\
\label{ccc11}&:=& T_1(n^\star) + T_2(n^\star) + T_3(n^\star).
\end{eqnarray}
Assumption~\textbf{(M3$^{\prime}$)} and the fact that $\sup_{\ell\ge 1}|\eta^p(\ell)|\le 1$ for every $p\in \mathcal{J}$ yield $\sum_{\ell=1}^\infty\alpha^{q_1}_{q,q_1}(\ell)\eta^{q_1}(\ell) < \infty$ so that $T_3(n^\star)\rightarrow 0$ in probability as $n^\star\rightarrow \infty$. Also from~\textbf{(M3$^{\prime}$)} we have
$\{\sum_{\ell=1}^\infty |U^{n_{q}}_{n_{q_1}}(n^\star,\ell)|^2(V^{n_{q_1}}(\ell))^{-1}; n^\star\in \mathbb{N}^J\}$ is bounded in probability and therefore Remark~\ref{hestimator} implies

$$|T_1(n^\star)|\le \sup_{\ell\ge 1}|\hat{\eta}^{n_{q_1}}(\ell)-\eta^{q_1}(\ell)|\sum_{\ell=1}^\infty
\frac{|U^{n_{q}}_{n_{q_1}}(n^\star,\ell)|^2}{V^{n_{q_1}}(\ell)}\rightarrow 0$$
in probability as $n^\star\rightarrow \infty$. The assertion that $T_2(n^\star)\rightarrow 0$ is a direct consequence of~\textbf{(M3$^{\prime}$)}. By using exactly the same above argument for the other terms in the difference $\sum_{\ell=1}^{d^u_{n^\star}}\hat{\phi}^2_{q,n^\star}(\ell)-\sum_{\ell=1}^\infty\phi^2_q(\ell)$ and $\sum_{\ell=1}^{d^u_{n^\star}}\hat{\psi}_{n^\star}(k,r,\ell)- \sum_{\ell=1}^\infty\psi(k,r,\ell)$ together with the correspondent assumptions~\textbf{(M4$^{\prime}$-H1$^{\prime}$)}, we arrive at~(\ref{ccc1}) and~(\ref{ccc2}). This allows us to conclude the proof.

\section{Proof of Theorem~\ref{main_asymptotic_result}}\label{proofTH2}

We start with the simplest case $d^u< \infty$ and to shorten notation we denote $W^{n^\star}_q(i)$ the $q$-th coordinate of the vector $\tilde{W}^{n^\star}(i)$ introduced in the proof of Theorem~\ref{Conv_LR_infinity}. Throughout this proof, any element $ x\in\mathbb{R}^p$ for $p< \infty$ is identified as an element of $\ell^2$ in the obvious way, $C$ is a constant which may defer from line to line and we set $k(p,q)=q-p+1$ for any $1\le p\le q < \infty$. We can write
$$GET(n^\star, \phi, p,q)=\big(M_{p}\tilde{W}^{n^\star}(p),\ldots, M_{q}\tilde{W}^{n^\star}(q)\big)$$
where $M_j:=diag~(\phi_1(j),\ldots, \phi_J(j));~j\ge 1$. A direct application of Theorem~\ref{Conv_LR_infinity} yields $M_i\tilde{W}^{n^\star}(i)\rightarrow N(0,M_i\Gamma(i)M_i)$ weakly as $n^\star\rightarrow \infty$ for each $i\ge d^l$ and more importantly, for every $a\in \mathbb{R}^{Jk(d^l,i)}$ we have

$$\big\langle GET(n^\star,\phi,d^l,i), a\big\rangle_{\mathbb{R}^{Jk(d^l,i)}}\rightarrow N\Big(0,\langle \mathcal{Y}(d^l,i)a,a\rangle_{\mathbb{R}^{Jk(d^l,i)}} \Big)$$
weakly as $n^\star\rightarrow \infty$. This shows that for each $i\ge d^l$

\begin{equation}\label{fconvv1}
GET(n^\star,\phi,d^l,i)\rightarrow N(0,\mathcal{Y}(d^l,i))
\end{equation}
and
\begin{equation}\label{fconvv}
\|GET(n^\star,\phi,d^l,i)\|^2_{\ell^2}\rightarrow \sum_{s=1}^{k(d^l,i)}\sum_{q=1}^J\lambda_{sq}\chi^2_{sq}
\end{equation}
weakly as $n^\star\rightarrow \infty$, where $\{\lambda_{sq} ; 1\le s\le k(d^l,i), q=1, \ldots , J\}$ are the eigenvalues of $\mathcal{Y}(d^l,i)$ and $\{\chi_{sq}^2; 1\le s\le k(d^l,i), q=1, \ldots , J\}$ is an i.i.d subset of chi-squared random variables with one degree of freedom. Next, we have to play with the estimators $\hat{\phi}_{n^\star},$ $d^l_{n^\star}$ and $d^u_{n^\star}$. Triangle inequality yields

\[
\|GET(n^\star, \phi,d^l,d^u)  - GET(n^\star, \hat{\phi}_{n^\star},d^l_{n^{\star}},d^u_{n^\star})\|_{\ell^2} \leq
\]

\[
\|GET(n^\star, \hat{\phi}_{n^\star},d^l,d^u)  - GET(n^\star, \phi,d^l,d^u)\|_{\ell^2} +
  \]

\[
\|GET(n^\star, \hat{\phi}_{n^\star},d^l,d^u)  - GET(n^\star, \hat{\phi}_{n^\star},d^l_{n^{\star}},d^u_{n^{\star}})\|_{\ell^2}=:
\]
\[
T_1(n^\star) + T_2(n^\star).
\]
For a given $q\in\mathcal{J}$ and $n^\star\in\mathbb{N}^J$, $\{W^{n^\star}_q; i\ge 1\}$ is an $\mathbb{F}$-martingale and therefore Doob's maximal inequality and assumptions \textbf{(M3$^{\prime}$,M4$^{\prime}$)} yield

\[
 \mathbb{E} \sup_{i\ge 1}|W^{n^\star}_q(i)|^2 \le C \mathbb{E} \big\langle W^{n^\star}_q \rangle (\infty) \leq\]

 \[
 \mathbb{E} \sum_{\ell = 1}^{\infty}\sum_{q_1\neq q}  \frac{|U^{n_q}_{n_{q_1}}(n^\star,\ell)|^2}{V^{n_{q_1}}(\ell)} +
\mathbb{E} \sum_{\ell = 1}^{\infty}\sum_{q_1\neq q}  \frac{|U^{n_q}_{n_{q_1}}(n^\star,\ell)|^2}{V^{n_{q}}(\ell)} +
\]

\begin{equation}\label{jh}
2 \mathbb{E} \sum_{\ell = 1}^{\infty}\sum_{(q_1,q_2) \in A_q} \frac{|U^{n_q}_{n_{q_1}}(n^\star,\ell)U^{n_q}_{n_{q_2}}(n^\star,\ell)|}{V^{n_q}(\ell)} \le C;~\forall n^\star\in \mathbb{N}^J.
\end{equation}
Estimate~(\ref{jh}) implies in particular that $\{\sup_{i\ge 1}|W^{n^\star}_q(i)|^2;n^\star\in\mathbb{N}^J\}$ is bounded in probability for each $q\in\mathcal{J}$. Therefore,  Lemma~\ref{con_ell_2} yields

\begin{equation}\label{get1}
T_1(n^\star)^2 \le \sum_{s=1}^{k(d^l,d^u)}\sum_{q=1}^J|\hat{\phi}_{q,n^\star}(s)-\phi_q(s)|^2\sum_{q=1}^J\sup_{i\ge 1}|W^{n^\star}_q(i)|^2\rightarrow 0
\end{equation}
in probability as $n^\star\rightarrow \infty$. By the very definition,

\begin{equation}\label{get2}
 T^2_2(n^\star)
:=\left\{
\begin{array}{crcr}
\sum_{\ell=k(d^l_{n^\star}, d^u_{n^\star})}^{k(d^l,d^u)}\sum_{q=1}^J |\hat{\phi}_{q,n^\star}(\ell)W^{n^\star}_q(\ell)|^2~&;&k(d^l_{n^\star}, d^u_{n^\star}) \le k(d^l,d^u)\\
\sum_{\ell=k(d^l,d^u)}^{k(d^l_{n^\star}, d^u_{n^\star})}\sum_{q=1}^J |\hat{\phi}_{q,n^\star}(\ell)W^{n^\star}_q(\ell)|^2~&;&k(d^l_{n^\star}, d^u_{n^\star}) > k(d^l,d^u).
\end{array}
\right.
\end{equation}
So $T_2(n^\star)\rightarrow 0$ in probability as $n^\star\rightarrow \infty$. Summing up the above estimates we  conclude~(\ref{etlimit}) and~(\ref{astat}) when $d^u< \infty$. Let us now treat the case $d^u=\infty$. At first, we notice that $\mathcal{Y}(d^l,d^u)$ is a nuclear operator. For a given $i\ge d^l$, let $\{\lambda_{sq} ; 1\le s\le k(d^l,i), q=1, \ldots , J\}$ be the eigenvalues of $\mathcal{Y}(d^l,i)$ and let $\{\chi_{sq}^2; 1\le s\le k(d^l,i), q=1, \ldots , J\}$ be an i.i.d subset of chi-squared random variables with one degree of freedom. Convergence~(\ref{fconvv}), properties~\textbf{(M3$^{\prime}$-M4$^{\prime}$)} and estimate~(\ref{jh}) yield

\begin{equation}\label{trace}
\sum_{s=1}^{Jk(d^l,i)}\lambda_{s} = \mathbb{E} \sum_{s=1}^{Jk(d^l,i)} \lambda_{s} \chi_{s}^2 \leq \liminf_{n^{\star}\rightarrow \infty} \mathbb{E} \|GET(n^\star,\phi, d^l,i)\|^2_{\ell^2} \le C\sum_{\ell\ge d^l}^\infty\sum_{q=1}^{J}\phi^2_q(\ell),
\end{equation}
for every $i\ge d^l$ and hence~(\ref{trace}) yields $Tr~(\mathcal{Y}(d^l,d^u))=\sum_{s=1}^\infty\lambda_s < \infty$. Let us now check tightness of the family $\{GET(n^\star,d^l,d^u);n^\star\in\mathbb{N}^J \}$. From~(\ref{jh}) and~\textbf{(M3$^{\prime}$-M4$^{\prime}$)}, we have

\begin{equation}\label{tight1}
\mathbb{E}\|GET(n^\star,d^l,d^u)\|^2_{\ell^2}\le \sum_{q=1}^J\mathbb{E}\sup_{s\ge 1}|W^{n^\star}_q(s)|^2\times\sum_{s=1}^\infty\sum_{q=1}^J|\phi_q(s)|^2\le C
\end{equation}
for every $n^\star\in \mathbb{N}^J$. In particular,~(\ref{jh}) yields

\begin{eqnarray}
\label{tight2}
\sup_{n^\star\in\mathbb{N}^J}\mathbb{E}\sum_{s=N}^\infty\sum_{q=1}^J|\phi_q(s)W^{n^\star}_q(s)|^2&\le& \sum_{q=1}^J\mathbb{E}\sup_{s\ge 1}|W^{n^\star}_q(s)|^2\times\sum_{r=N}^\infty\sum_{q=1}^J|\phi_q(r)|^2\\
\nonumber& &\\
\nonumber&\le& C\sum_{q=1}^J\sum_{r=N}^\infty|\phi_q(r)|^2\rightarrow 0,
\end{eqnarray}
as $N\rightarrow \infty$. Hence,~\textbf{(M3$^{\prime}$-M4$^{\prime}$)},~(\ref{tight1}) and~(\ref{tight2}) allow us to conclude the relatively weak compactness which together with the weak convergence of the finite-dimensional projections~(\ref{fconvv1}) imply
\begin{equation}\label{almost1}
GET(n^\star, \phi,d^l,d^u)\rightarrow N(0,\mathcal{Y}(d^l,d^u))
\end{equation}
weakly as $n^\star\rightarrow \infty$. It remains to play with the estimators, but this is a straightforward consequence of the previous arguments. In fact, triangle inequality yields

\[
\big|\|GET(n^\star, \phi,d^l,d^u)\|^2_{\ell^2}  - \|GET(n^\star, \hat{\phi}_{n^\star},d^l_{n^{\star}},d^u_{n^\star})\|^2_{\ell^2} \big|\leq
\]
\[
\big|\|GET(n^\star, \hat{\phi}_{n^\star},d^l,d^u)\|^2_{\ell^2}  - \|GET(n^\star, \phi,d^l,d^u)\|^2_{\ell^2}\big| +
\]

\[
\big|\|GET(n^\star, \hat{\phi}_{n^\star},d^l,d^u)\|^2_{\ell^2}  - \|GET(n^\star, \hat{\phi}_{n^\star},d^l_{n^{\star}},d^u_{n^{\star}})\|^2_{\ell^2}\big|=:
\]
\[
T_3(n^\star) + T_4(n^\star).
\]
The same arguments given in~(\ref{ccc1}) and~(\ref{ccc11}) allows us to get $\sum_{s=1}^{\infty}\sum_{q=1}^J|\hat{\phi}^2_{q,n^\star}(s)-\phi^2_q(s)|\rightarrow 0$ in probability as $n^\star\rightarrow \infty$ and since $\sup_{i\ge 1}|W^{n^\star}_q(i)|^2$ is bounded in probability, we can safely conclude
\begin{equation}\label{get33}
T_3(n^\star) \le \sum_{s=1}^{\infty}\sum_{q=1}^J|\hat{\phi}^2_{q,n^\star}(s)-\phi^2_q(s)|\sum_{q=1}^J\sup_{i\ge 1}|W^{n^\star}_q(i)|^2\rightarrow 0
\end{equation}
in probability as $n^\star\rightarrow \infty$. By the same reason,

$$T_4(n^\star)\le \sum_{s=k(d^l_{n^\star},d^u_{n^\star})+1}^\infty\sum_{q=1}^J |\hat{\phi}_{q,n^\star}(s)
W^{n^\star}_q(s)|^2\rightarrow 0$$
in probability as $n^\star\rightarrow \infty$. Convergence~(\ref{almost1}) jointly with $T_3(n^\star)+T_4(n^\star)\rightarrow 0$ in probability as $n^\star\rightarrow \infty$ allow us to conclude the proof.

\section{Proof of Proposition~\ref{HAass}}\label{proofProp}
At first, we check (\textbf{M1-M4}) and (\textbf{H1}). Let us fix a positive integer $\ell\ge 1$. Assumption \textbf{(M1)} is obvious. If $\delta > 0$ then we shall use the sample growth condition and the definition of $u$ to get the following estimate

\begin{eqnarray*}
n_q \Bigg|\frac{U_{n_{q_1}}^{n_q}(n^\star,i)}{V^{n_q}(i)}\Bigg|^{2+\delta}&\le& \Big(\frac{n_q}{n}\Big) n^{-\delta / 2}|u(n^\star,\ell)|^{2+ \delta}\\
& &\\
&\rightarrow& 0\quad \text{as}~n^\star\rightarrow \infty,
\end{eqnarray*}

\noindent for any $q\ne q_1$ in $\mathcal{J}$. This shows that assumption \textbf{(M2)} is satisfied. For a given $q\neq q_1$ in $\mathcal{J}$ and $q_2 = q$ we have

\begin{equation}\label{cm3}
\frac{|U^{n_{q}}_{n_{q_1}}(n^\star,\ell)|^2}{V^{n_{q_2}}(\ell)}=\Bigg|u(n^ \star,\ell)\frac{V^{n_{q_1}}(\ell)}{V^{n^\star(\ell)}}\Bigg|^2 \frac{n_q}{n} \frac{V^{n_q}(\ell)}{n_q}.
\end{equation}
For $q_2=q_1$, we have a similar expression. Now, if $q_2 \in \mathcal{J}$~$(q_2 \neq q, q_2 \neq q_1)$ then we shall write

\begin{equation} \label{cm4}
\frac{U^{n_{q}}_{n_{q_1}}(n^\star,\ell)U^{n_q}_{n_{q_2}}(n^\star,\ell)}{V^{n_q}(\ell)}
= |u(n^\star,\ell)|^2\frac{ V^{n_{q_2}}(\ell)}{V^{n^\star}(\ell)} \frac{V^{n_{q_1}}(\ell)}{V^{n^\star}(\ell)}  \frac{n_q}{n}  \frac{V^{n_q}(\ell)}{n_{q}}
\end{equation}

Identities~(\ref{cm3}) and~(\ref{cm4}) allow us to use again the sample growth condition, the definition of $u$ and the binomial property to get assumption (\textbf{M3}) and (\textbf{M4}). In fact,

$$
\alpha^{q_2}_{q,q_1}(\ell) =  \left( \omega (\ell) \frac{b_{q_1} \theta^{q_1} (\ell)}{\sum_{p=1}^J b_p \theta^p (\ell)} \right)^2 b_q \theta^q (\ell)
$$ and

$$
\beta^{q}_{q_1,q_2}(\ell) =  \omega^2 (\ell) \frac{b_{q_2} \theta^{q_2} (\ell)}{\sum_{p=1}^J b_p \theta^p (\ell)} \frac{b_{q_1} \theta^{q_1} (\ell)}{\sum_{p=1}^J b_p \theta^p (\ell)} b_q \theta^q (\ell).
$$

Let us now check \textbf{(H1)}. If $q_1\neq k \in \mathcal{J}$, then we shall write

\begin{equation} \label{cm5}
\frac{U^{n_{k}}_{n_{q_1}}(n^\star,\ell)U^{n_r}_{n_{q_1}}(n^\star,\ell)}{V^{n_{q_1}}(\ell)}
= |u(n^\star,\ell)|^2\frac{V^{n_k}(\ell)}{V^{n^\star}(\ell)} \frac{V^{n_{r}}(\ell)}{V^{n^\star}(\ell)} \frac{n_{q_1}}{n}  \frac{V^{n_{q_1}}(\ell)}{n_{q_1}}.
\end{equation}

\noindent For a given $r\neq k$ in $\mathcal{J}$ and $q_1\neq k$ we shall write

\begin{equation} \label{cm6}
\frac{U^{n_{k}}_{n_{q_1}}(n^\star,\ell)U^{n_r}_{n_{k}}(n^\star,\ell)}{V^{n_{k}}(\ell)}
= |u(n^\star,\ell)|^2\frac{V^{n_r}(\ell)}{V^{n^\star}(\ell)} \frac{V^{n_{k}}(\ell)}{V^{n^\star}(\ell)} \frac{n_{q_1}}{n}  \frac{V^{n_{q_1}}(\ell)}{n_{q_1}}.
\end{equation}


\noindent Identities~(\ref{cm5}) and~(\ref{cm6}) allow us to use again the sample growth condition, the binomial property and the definition of $u$ to get assumption \textbf{(H1)}. In this case,

\begin{equation} \label{covariance_neg}
\eta^{k,r}_{q_1}(\ell) = \gamma^{k,r}_{q_1}(\ell) =  \omega^2 (\ell) \frac{b_{k} \theta^{k} (\ell)}{\sum_{p=1}^J b_p \theta^p (\ell)} \frac{b_{r} \theta^{r} (\ell)}{\sum_{p=1}^J b_p \theta^p (\ell)} b_{q_1} \theta^{q_1} (\ell).
\end{equation}
It remains to check assumptions \textbf{(M3$^{\prime}$-M4$^{\prime}$)} and \textbf{(H1$^{\prime}$)}. For the assumption \textbf{(M3$^{\prime}$)}, we notice that if we take $q_2 =q$, there exists a positive constant $C$ such that

%

\begin{equation} \label{fiiu}
 \frac{|U^{n_{q}}_{n_{q_1}}(n^\star,\ell)|^2}{V^{n_{q}}(\ell)} \leq C  \frac{V^{n_q}(\ell)}{n_q}, \quad \ell \geq 1.
\end{equation}
Therefore, we have that

$$ \sum_{\ell=1}^{\infty} \limsup_{n^\star}\mathbb{E}\frac{|U^{n_{q}}_{n_{q_1}}(n^\star,\ell)|^2}{V^{n_{q}}(\ell)} \leq C \sum_{\ell=1}^{\infty} \limsup_{n^\star}\mathbb{E}  \frac{V^{n_q}(\ell)}{n_q} \leq C \sum_{\ell=1}^{\infty} \theta^q (\ell) < \infty$$

Obviously,

\begin{equation} \label{ineq_weigthed_proc} \Bigg|\frac{|U^{n_{q}}_{n_{q_1}}(n^\star,\ell)|^2}{V^{n_{q_2}}(\ell)}-\alpha^{q_2}_{q,q_1}(\ell)\Bigg| \leq \frac{|U^{n_{q}}_{n_{q_1}}(n^\star,\ell)|^2}{V^{n_{q_2}}(\ell)} +  \alpha^{q_2}_{q,q_1}(\ell).
\end{equation}
and from~(\ref{cm3}) and the binomial property we actually have
\begin{equation}\label{cm7}
\lim_{n^\star\rightarrow \infty}\mathbb{E}\Bigg|\frac{|U^{n_{q}}_{n_{q_1}}(n^\star,\ell)|^2}{V^{n_{q_2}}(\ell)}
-\alpha^{q_2}_{q,q_1}(\ell)\Bigg|=0
\end{equation}
for each $\ell\ge 1$. For a given $\varepsilon > 0$

\begin{eqnarray}
\nonumber\lim_{n^{\star} \rightarrow \infty} \mathbb{P} \left[ \sum_{\ell=1}^{\infty}\Bigg|\frac{|U^{n_{q}}_{n_{q_1}}(n^\star,\ell)|^2}{V^{n_{q_2}}(\ell)}-\alpha^{q_2}_{q,q_1}(\ell)\Bigg| > \varepsilon \right] &\leq& \lim_{n^{\star} \rightarrow \infty} \frac{1}{\varepsilon} \sum_{\ell=1}^{\infty}\mathbb{E}\Bigg|\frac{|U^{n_{q}}_{n_{q_1}}(n^\star,\ell)|^2}{V^{n_{q_2}}(\ell)}-
\alpha^{q_2}_{q,q_1}(\ell)\Bigg|\\
\nonumber& &\\
\label{cm8}&=&\frac{1}{\varepsilon}\sum_{\ell=1}^{\infty} \lim_{n^{\star}\rightarrow \infty}\mathbb{E} \Bigg|\frac{|U^{n_{q}}_{n_{q_1}}(n^\star,\ell)|^2}{V^{n_{q_2}}(\ell)}-\alpha^{q_2}_{q,q_1}(\ell)\Bigg| \\
\nonumber& &\\
\nonumber&=& 0.
\end{eqnarray}
The justification of the limit into the series in~(\ref{cm8}) is due to~(\ref{fiiu}) and (\ref{cm7}) which gives a constant $C>0$ such that

$$\mathbb{E}\frac{|U^{n_{q}}_{n_{q_1}}(n^\star,\ell)|^2}{V^{n_{q_2}}(\ell)}\le C\mathbb{E}\frac{V^{n_{q}}(\ell)}{n_{q}}=C\theta^{q}(\ell); \ell \ge 1,
$$ and

$$\alpha^{q_2}_{q,q_1}(\ell) = \lim_{ n^{\star} \rightarrow \infty } \frac{|U^{n_{q}}_{n_{q_1}}(n^\star,\ell)|^2}{V^{n_{q_2}}(\ell)} \le  C\theta^{q}(\ell); \ell \ge 1,
$$
for every $n^\star\in\mathbb{N}^J$ where $\theta^q(\cdot)\in \ell^1(\mathbb{N})$ (by assumption $X^q$ is integrable). We can apply the same arguments to check assumptions \textbf{(M4$^{\prime}$)} and \textbf{(H1$^{\prime}$)}. This concludes the proof.

\section{Proof of Theorem~\ref{main_CVM}}\label{proofTH3}
The proof is an almost direct consequence of Theorem~\ref{main_asymptotic_result}. Convergence~(\ref{BCVM}) is consequence of Theorem~\ref{main_asymptotic_result} and the only statement which has to be detailed is convergence~(\ref{apr_CVM}) when $d^u=\infty$. We take $d^u_{n^\star}$ given by~~(\ref{dustar}) and without any loss of generality and to simplify notation, we set $d^l_{n^\star}=1$. The arguments for general $d^l_{n^\star}$ follow easily from this case. Let us define

\[
A(n^{\star}, N) := \left\{\hat{\mathcal{Y}}_0( 1,N \wedge d_{n^{\star}}^u)~\text{is non-negative}\right\}
\]
and we notice that $A(n^{\star})=  \cap_{N=1}^{\infty} A(n^{\star}, N)$. We claim that $\lim_{n^\star\rightarrow \infty} \mathbb{P}(A(n^{\star})) = 1$. In fact, a basic inequality among compact operators~(see e.g~\cite{gohberg}) yields the following a.s estimate

\begin{equation}\label{eigin}
|\lambda_{sq} - \hat{\lambda}_{sq}|\le \|\hat{\mathcal{Y}}_0(1,d^u_{n^\star}\wedge N) - \mathcal{Y}_0(1,d^u_{n^\star}\wedge N)\|
\end{equation}
for every $q=1,\ldots, J-1$ and $1\le s\le d^u_{n^\star}\wedge N$. Here $\|\cdot\|$ stands the strong norm over the space of bounded operators in $\ell^2$. Since $\inf\{\lambda_{sq}; 1\le s\le N, 1\le q\le J-1 \}\ge 0$, we may use~(\ref{eigin}) together with Lemma~\ref{con_ell_2} and Theorem~\ref{Conv_LR_infinity} to conclude that

$$\lim_{n^\star\rightarrow \infty}\mathbb{P}\Big(\bigcap_{N=1}^kA(n^\star,N) \Big)=1~\text{for each}~k\ge 1.$$
Hence, we do have

\begin{eqnarray}
\label{lim A}\lim_{n^\star\rightarrow \infty}\mathbb{P}\Big( A(n^\star)\Big)&=&\lim_{n^\star\rightarrow \infty}\lim_{k\rightarrow \infty}\mathbb{P}\Big(\bigcap_{N=1}^kA(n^\star,N)\Big)\\
\nonumber& &\\
\nonumber&=&\lim_{N\rightarrow \infty}\lim_{n^\star\rightarrow \infty}\mathbb{P}\Big(\bigcap_{N=1}^kA(n^\star,N)\Big)=1.
\end{eqnarray}

For a given i.i.d sequence $\{X_{i}\}_{i=1}^\infty$ of real-valued Gaussian variables $N(0,1)$ and positive integer $N\ge 1$, let us define the following sequences of $\ell^2$-valued random variables

$$\mathcal{X}(sq):=\sqrt{\lambda}_{sq}X_{sq}; q=1,\ldots, J-1, s\ge 1;$$

$$
\mathcal{X}^N(sq) :=\left\{
\begin{array}{crcr}
\sqrt{\lambda}_{sq}X_{sq}~&;&q=1,\ldots, J-1; 1\le s\le N \\
0~&;& s> N;
\end{array}
\right.
$$

$$
\hat{\mathcal{X}}_{n^\star}(sq) :=\left\{
\begin{array}{crcr}
\sqrt{\hat{\lambda}}_{sq}X_{sq} 1\!\!1_{A(n^{\star})}~&;&q=1,\ldots, J-1; 1\le s\le d^u_{n^\star} \\
0~&;&~s > d^u_{n^\star};
\end{array}
\right.
$$

$$
\hat{\mathcal{X}}^N_{n^\star}(sq) :=\left\{
\begin{array}{crcr}
\sqrt{\hat{\lambda}}_{sq}X_{sq} 1\!\!1_{A(n^{\star})}~&;&q=1,\ldots, J-1; 1\le s\le d^u_{n^\star}\wedge N \\
0~&;& s > d^u_{n^\star}\wedge N.
\end{array}
\right.
$$
Convergence~(\ref{lim A}) and the inequality~(\ref{eigin}) yield

$$
\|\hat{\mathcal{X}}^N_{n^\star}-\mathcal{X}^N\|^2_{\ell^2}\rightarrow 0
$$
in probability as $n^\star\rightarrow \infty$ for each $N\ge 1$. Of course, $\lim_{N\rightarrow \infty}\mathcal{X}^N=\mathcal{X}$ in probability in $\ell^2$. The strategy is to prove that for a given $\varepsilon > 0$

\begin{equation}\label{vvvv2}
\lim_{N\rightarrow \infty}\limsup_{n^\star\rightarrow \infty}\mathbb{P} \big\{\|\hat{\mathcal{X}}^N_{n^\star} - \hat{\mathcal{X}}_{n^\star} \|_{\ell^2} > \varepsilon \big\}=0.
\end{equation}
Under~(\ref{vvvv2}), we may exchange the iterated weak $\ell^2$-limits $\lim_{N}\lim_{n^\star}\hat{\mathcal{X}}^N_{n^\star} = \lim_{n^\star}\lim_{N}\hat{\mathcal{X}}^N_{n^\star}$ which allow us to conclude~(\ref{apr_CVM}), i.e., $\|\hat{\mathcal{X}}_{n^\star}\|^2_{\ell^2}\rightarrow \|\mathcal{X}\|^2_{\ell^2}$ weakly as $n^\star\rightarrow \infty$.
By the very definition,

\begin{equation}\label{ff1}
\|\hat{\mathcal{X}}^N_{n^\star} - \hat{\mathcal{X}}_{n^\star} \|^2_{\ell^2}=\sum_{q=1}^{J-1}\sum_{s=N+ 1}^{d^u_{n^\star}}\hat{\lambda}_{sq}\chi^2_{sq} 1\!\!1_{A(n^{\star})} ;~1\le N < d^u_{n^\star},
\end{equation}
and $\hat{Y}_0(1, d^u_{n^\star})1\!\!1_{A(n^{\star})}$ is non-negative a.s for every $n^\star\in \mathbb{N}^J$. In the sequel, we assume that $\hat{\lambda}$ are enumerated with algebraic multiplicities taked into account. By the Lidskii trace theorem

$$\sum_{s=1}^{d^u_{n^\star}}\sum_{q=1}^{J-1}\hat{\lambda}_{sq} 1\!\!1_{A(n^{\star})} = Tr~\hat{Y}_0(1,d^u_{n^\star}) 1\!\!1_{A(n^{\star})}$$
a.s for each $n^\star\in \mathbb{N}^J$. In particular, from the fact that $\hat{Y}_0(1, d^u_{n^\star})1\!\!1_{A(n^{\star})}$ is non-negative a.s we get the following bound [see e.g~Corollary 3.7.p.56 in~\cite{gohberg}]

\begin{eqnarray}
\nonumber\sum_{s=N+1}^{d^u_{n^\star}}\sum_{q=1}^{J-1}\hat{\lambda}_{sq} 1\!\!1_{A(n^{\star})} &\le& \sum_{j=N+1}^{d^u_{n^\star}}Tr~\hat{M}_0(n^\star,j)\hat{\Gamma}_0(n^\star,j) \hat{M}_0(n^\star,j)1\!\!1_{A(n^{\star})} \\
\nonumber& &\\
\label{lid}&=&\sum_{j=N+1}^{d^u_{n^\star}}\sum_{q=1}^{J-1}\Bigg\{\sum_{\ell=1}^{j}\hat{\phi}^2_{q,n^\star}(\ell)\Bigg\}\times \hat{\phi}^2_{q,n^\star}(j)1\!\!1_{A(n^{\star})}~a.s,
\end{eqnarray}
for $n^\star\in \mathbb{N}^J$. We know there exists $C>0$ (which only depends on $(b_p)_{p=1}^J$, see~(\ref{fiiu})) such that

\begin{eqnarray}
\nonumber\sum_{j=N+1}^{d^u_{n^\star}}\sum_{q=1}^{J-1}\Bigg\{\sum_{\ell=1}^{j}\hat{\phi}^2_{q,n^\star}(\ell)\Bigg\}
\times\hat{\phi}^2_{q,n^\star}(j)\chi^2_{jq}1\!\!1_{A(n^{\star})} &\leq& C \sum_{j=N+1}^{d^u_{n^\star}}\sum_{q=1}^{J-1}\Bigg\{\sum_{\ell=1}^{j}\frac{V^{n_q}(\ell)}{n_q}\Bigg\}\\
\nonumber& &\\
\label{RER1}&\times&\frac{V^{n_q}(j)}{n_q}\chi^2_{jq}1\!\!1_{A(n^{\star})}
\end{eqnarray}
a.s for every~$n^\star\in \mathbb{N}^J$. In particular, we may assume that $\{\chi^2_{jp};j\ge N+1, 1\le p\le J-1\}$ are independent of $V^{n_q}$ for every $q\in \mathcal{J}$ to get

\begin{equation}\label{RER2}
\mathbb{E} \sum_{j=N+1}^{d^u_{n^\star}}\sum_{q=1}^{J-1}\Bigg\{\sum_{\ell=1}^{j}\frac{V^{n_q}(\ell)}{n_q}\Bigg\}
\times\frac{V^{n_q}(j)}{n_q}\chi^2_{jq} \leq \sum_{q=1}^{J-1} \sum_{j=N+1}^{\infty} \sum_{\ell=1}^{j} \mathbb{E} \frac{V^{n_q}(\ell)}{n_q}
\times\frac{V^{n_q}(j)}{n_q}.
\end{equation}
By taking advantage of the independence of the samples, a straightforward calculation yields

\begin{eqnarray}
\label{indepcov}\mathbb{E} V^{n_q}(\ell) V^{n_q}(j) &=&n_q \theta^q(j) \left( 1 -  \theta^q(\ell)  \right) + n_q \theta^q(\ell) n_q \theta^q(j)\\
\nonumber& &\\
\nonumber&\le& n_q \theta^q(j) +  n_q \theta^q(\ell) n_q \theta^q(j)
\end{eqnarray}
for every $n_q\ge 1,~q\in\mathcal{J}$ and $1\le \ell \le j < \infty$. Hence, from~(\ref{indepcov}) we obtain the following bound

$$
\sum_{q=1}^{J-1} \sum_{j=N+1}^{\infty} \sum_{\ell=1}^{j} \mathbb{E} \frac{V^{n_q}(\ell)}{n_q}
\times\frac{V^{n_q}(j)}{n_q}  \leq \sum_{q=1}^{J-1} \sum_{j=N+1}^{\infty} \sum_{\ell=1}^{j} \left[\frac{\theta^q(j)}{n_q} +  \theta^q(\ell) \theta^q(j) \right] =
$$

$$
\sum_{q=1}^{J-1} \sum_{j=N+1}^{\infty} j \frac{\theta^q(j)}{n_q}  + \sum_{q=1}^{J-1} \sum_{j=N+1}^{\infty} \sum_{\ell=1}^{j} \left[ \theta^q(\ell) \theta^q(j) \right] \leq
$$

\begin{equation}\label{RER3}
\sum_{q=1}^{J-1} \frac{1}{n_q}\sum_{j=N+1}^\infty  j\theta^q(j)  + \sum_{q=1}^{J-1} \sum_{j=N+1}^{\infty} \theta^q(j)\mathbb{E}X^q,
\end{equation}
where $\sum_{j=N+1}^\infty  j\theta^q(j) < \infty$ for every $N\ge 1$ due to the integrability assumption $X^q \in L^2(\mathbb{P});~q=1,\ldots, J-1$. Summing up the above steps~(\ref{RER1}),~(\ref{RER2}) and~(\ref{RER3}), for a given $\varepsilon > 0$ the following estimate holds

\begin{eqnarray}
\nonumber\limsup_{n^\star\rightarrow \infty}\mathbb{P}\big\{\|\hat{\mathcal{X}}^N_{n^\star}-\hat{\mathcal{X}}_{n^\star}\|^2_{\ell^2} > \varepsilon\big\}&\le&\limsup_{n^\star\rightarrow \infty}
\frac{1}{\varepsilon}\mathbb{E}\sum_{s=N+1}^{d^u_{n^\star}}\sum_{q=1}^{J-1}\hat{\lambda}_{sq}\chi^2_{sq}1\!\!1_{A(n^{\star})}\\
\nonumber& &\\
\nonumber&\le&\frac{1}{\varepsilon}\limsup_{n^\star\rightarrow \infty}\mathbb{E}\sum_{j=N+1}^{d^u_{n^\star}}
\nonumber\sum_{q=1}^{J-1}\Bigg\{\sum_{\ell=1}^{j}\hat{\phi}^2_{q,n^\star}(\ell)\Bigg\}\times \hat{\phi}^2_{q,n^\star}(j)\chi^2_{jq}\\
\nonumber&\le&\frac{C}{\varepsilon} \limsup_{n^\star\rightarrow \infty}\sum_{q=1}^{J-1} \frac{1}{n_q}\sum_{j=N+1}^\infty  j\theta^q(j)\\
\nonumber& &\\
\label{lastestt}&+&\sum_{q=1}^{J-1} \sum_{j=N+1}^{\infty} \theta^q(j)\mathbb{E}X^q=\sum_{q=1}^{J-1} \sum_{j=N+1}^{\infty} \theta^q(j)\mathbb{E}X^q.
\end{eqnarray}
Finally, from~(\ref{lastestt}) and the integrability of $X^q$ we may conclude~(\ref{vvvv2}). This shows $\|\hat{\mathcal{X}}_{n^ \star}\|^2_{\ell^2}\rightarrow \|\mathcal{X}\|^2_{\ell^2}$ weakly as $n^\star\rightarrow \infty$ which allows us to conclude the proof.

\section{Appendix}\label{appendixlogrank}
In this appendix, we provide the proofs of Lemmas~\ref{multimardiff},~\ref{con_ell_2},~\ref{mainres1} and Proposition~\ref{multlema}.

\subsection{Proof of Lemma~\ref{multimardiff}}
\begin{proof}
Let us fix $\ell\ge 1$, $n^\star\in \mathbb{N}^J$ and $k\in \mathcal{J}$. By construction, we notice that

\begin{eqnarray}
\nonumber\mathbb{E}[\Delta R^k_m(\ell)|\mathcal{G}^{n^\star}_{m-1}(\ell)]&=&\mathbb{E}[\Delta R^k_m(\ell)|V^k_m(\ell)]\\
\nonumber & &\\
\nonumber &=& V^k_m(\ell)h^k(\ell)\\
\nonumber & &\\
\label{key1}&=& \mathbb{E}[\Delta N^k_m(\ell)|\mathcal{G}^{n^\star}_{m-1}(\ell)],\quad 1\le m\le n_k.
\end{eqnarray}

\noindent Therefore $\mathbb{E}[\Delta Y^k_m(\ell)|\mathcal{G}^{n^\star}_{m-1}(\ell)]=0$ for $1\le m\le n_k$. Assumption~\textbf{(M1)} allows us to conclude the proof.

\end{proof}

\subsection{Proof of Lemma~\ref{con_ell_2}}
\begin{proof}
At first, we notice that at the category $i$, the candidate variance $\phi^2_q(\ell)$ is well-defined for every $\ell\in \{1,\ldots, i \}$ and $q\in \mathcal{J}$ so let us fix such $\ell$ and a population $q$. In order to apply the classical martingale central limit theorem, we begin by verifying the Lindeberg condition. Indeed, it is sufficient to establish the conditional Liapunov condition

$$
\sum_{m=1}^{N_J} \mathbb{E} \left[  |\xi^{n^\star}_{m,q}(\ell)|^{2+ \delta} | \mathcal{G}^{n^\star}_{m-1} (\ell) \right] \rightarrow 0 \quad \mbox{in probability as}~n^\star\rightarrow \infty,
$$

\noindent for some $\delta>0$. By the very definition of $R^q_m$ and Doob-Meyer decomposition we have $|\Delta Y^{q}_{m}(\ell) |^{2+ \delta} \leq 2^{2+ \delta} $ a.s and hence we may use assumptions \textbf{(M1)} and \textbf{(M2)} to find a constant $C$ which only depends on $J$ and $\delta$ such that

\begin{eqnarray*}
\sum_{m=1}^{N_J}\mathbb{E}\big[|\xi^{n^\star}_{m,q}(\ell)|^{2+\delta}\big|\mathcal{G}^{n^\star}_{m-1}(\ell)\big] &\le& C \sum_{q_1\neq q}n_q\Bigg|\frac{U^{n_q}_{n_{q_1}}(n^\star,\ell)}{V^{n_q}(\ell)}\Bigg|^{2+\delta}\\
& &\\
&+& C \sum_{q_1\neq q}n_{q_1}\Bigg|\frac{U^{n_q}_{n_{q_1}}(n^\star,\ell)}{V^{n_{q_1}}(\ell)}\Bigg|^{2+\delta}\\
& &\\
&\rightarrow 0&~\text{as}~n^\star\rightarrow \infty.
\end{eqnarray*}

\noindent In order to shorten notation, let us define

\begin{equation}\label{pi}
\pi^q_{m,n^\star}(q_1,\ell):= \frac{U^{n_q}_{n_{q_1}}(n^\star,\ell)}{V^{n_q}(\ell)}\Delta Y^q_m(\ell); \quad q_1\neq q.
\end{equation}

\begin{equation}\label{lambda}
\lambda^q_{m,n^\star}(q_1,\ell):= \frac{U^{n_q}_{n_{q_1}}(n^\star,\ell)}{V^{n_{q_1}}(\ell)}\Delta Y^{q_1}_m(\ell); \quad q_1\neq q.
\end{equation}

The quadratic variation of the martingale $\sum_{m=1}^\cdot\xi^{n^\star}_{m,q}(\ell)$ at the point $N_J$ can be written as

\begin{equation}\label{imp1}
\sum_{m=1}^{N_J} \mathbb{E} \Big[|\xi^{n^\star}_{m,q}(\ell)|^{2} | \mathcal{G}^{n^\star}_{m-1} (\ell)\Big]=\end{equation}
$$\sum_{m=1}^{N_J}\sum_{q_1\neq q}\mathbb{E}\Big[|\pi_{m,n^\star}^q(q_1,\ell) - \lambda_{m,n^\star}^q(q_1,\ell)|^2|\mathcal{G}^{n^\star}_{m-1}(\ell)\Big]+$$
$$2\sum_{m=1}^{N_J}\sum_{(q_1,q_2) \in A_q}\mathbb{E}\Big[(\pi^q_{m,n^\star}(q_1,\ell) - \lambda^q_{m,n^\star}(q_1,\ell)) (\pi^q_{m,n^\star}(q_2,\ell) - \lambda^q_{m,n^\star}(q_2,\ell))\Big|\mathcal{G}^{n^\star}_{m-1}(\ell)\Big]$$
$$=:T_1(n^\star,\ell) + T_2(n^\star,\ell).$$

\noindent Assumption~\textbf{(M1)}, the independence of the random sample and~(\ref{key1}) yield

\begin{eqnarray}
\label{ast1}T_1(n^\star,\ell)&=&\sum_{q_1\neq q}\frac{|U_{n_{q_1}}^{n_q}(n^\star,\ell)|^2}{V^{n_{q_1}}(\ell)}h^{q_1}(\ell)[1-h^{q_1}(\ell)]\\
\nonumber& &\\
\nonumber&+&\sum_{q_1\neq q}\frac{|U_{n_{q_1}}^{n_q}(n^\star,\ell)|^2}{V^{n_q}(\ell)}h^{q}(\ell)[1-h^{q}(\ell)].
\end{eqnarray}

Assumption \textbf{(M3)} then yields

\begin{equation}\label{t1}
\lim_{n^\star\rightarrow \infty}T_1(n^\star,\ell)= \sum_{q_1\neq q}\alpha_{q,q_1}^{q_1}(\ell)h^{q_1}(\ell)[1-h^{q_1}(\ell)] + \sum_{q_1\neq q}\alpha_{q,q_1}^{q}(\ell)h^{q}(\ell)[1-h^{q}(\ell)]
\end{equation}

\noindent in probability. The arguments used for the first term can also be applied to $T_2(n^\star,\ell)$ and then we shall write

\begin{equation}\label{ast2}
T_2(n^\star,\ell)= 2\sum_{(q_1,q_2)\in A_q} \frac{U^{n_q}_{n_{q_1}}(n^\star,\ell)U^{n_q}_{n_{q_2}}(n^\star,\ell)}{V^{n_q}(\ell)}h^q(\ell)[1-h^q(\ell)].
\end{equation}

Assumption \textbf{(M4)} then yields

\begin{equation}\label{t2}
\lim_{n^\star\rightarrow \infty}T_2(n^\star,\ell)=2\sum_{(q_1,q_2)\in A_q} \beta^q_{q_1,q_2}(\ell)h^q(\ell)[1-h^q(\ell)].
\end{equation}

\noindent in probability. Summing up~(\ref{t1}) and~(\ref{t2}) we conclude that

\begin{equation}\label{d}
\sum_{m=1}^{N_J} \mathbb{E} \Big[|\xi^{n^\star}_{m,q}(\ell)|^{2} | \mathcal{G}^{n^\star}_{m-1} (\ell)\Big]\rightarrow \phi^2_q(\ell)~\text{in probability as}~n^\star\rightarrow \infty.
\end{equation}

Summing up the above steps, the martingale central limit theorem applied to $\{ \xi^{n^\star}_{m,q}(\ell); 1\le m \le N_J\}$ ensures the weak convergence~(\ref{w1}). The convergence $\lim_{n^\star\rightarrow\infty}\hat{\phi}^2_{q,n^\star}(\ell)=
\phi^2_q(\ell)$ in probability is a consequence of relations~(\ref{ast1}) and~(\ref{ast2}) combined with Remark~\ref{hestimator}.
\end{proof}

\subsection{Proof of Lemma~\ref{mainres1}}

\begin{proof}
We fix $q\in\mathcal{J}$ and a category $i$ such that $min_{1\le p\le J}\{\theta^p(i)\}>0$. Let us denote by $Z_q(\cdot)$ the weak limit of~(\ref{w1}) in Lemma~\ref{con_ell_2}

\begin{equation}\label{zq}
Z_q(\ell):=\lim_{n^\star\rightarrow \infty}\sum_{m=1}^{N_J}\xi^{n^\star}_{m,q}(\ell);\quad 1\le \ell \le i.
\end{equation}

\noindent We also set

\begin{equation}\label{wq}
W_q(i):=\sum_{\ell=1}^iZ_q(\ell).
\end{equation}

By the very definition

$$\sum_{\ell=1}^i\mathbb{E}[Z_q(\ell)|\mathcal{F}_{i-1}] = W_q(i-1) + \mathbb{E}[Z_q(i)|\mathcal{F}_{i-1}].$$

Since $\frac{U^{n_q}_{n_{q_1}}(n^\star,\cdot)}{V^{n_{q_2}}(\cdot)}$ is $\mathbb{F}$-predictable for any $q_1\neq q$ and $q_2\in \mathcal{J}$, we may use the martingale property of $\{Y^{n_{p}},~p\in \mathcal{J}\}$ and the definition of $\mathbb{F}$ to get

\begin{eqnarray}
\label{zmart}\mathbb{E}[Z_q(i)|\mathcal{F}_{i-1}]&=&\lim_{n^\star\rightarrow \infty}\sum_{q_1\neq q}\frac{U^{n_q}_{n_{q_1}}(n^\star,i)}{V^{n_{q}}(i)}\mathbb{E}[\Delta Y^{n_q}(i)|\mathcal{F}^q_{i-1}]\\
\nonumber& & \\
\nonumber&-&\lim_{n^\star\rightarrow \infty}\sum_{q_1\neq q}\frac{U^{n_q}_{n_{q_1}}(n^\star,i)}{V^{n_{q_1}}(i)}\mathbb{E}[\Delta
\nonumber Y^{n_{q_1}}(i)|\mathcal{F}^{q_1}_{i-1}]\\
\nonumber & &\\
\nonumber&=& 0\quad \text{(weak sense)}.
\end{eqnarray}

\noindent This shows that $W_q$ is an $\mathbb{F}$-martingale. At this point, from Lemma~\ref{con_ell_2} we only need to check that $W_q$ has independent increments. For this, we claim that $\langle W_q\rangle$ is a deterministic process. In order to shorten notation, let us define

$$\pi^q_{n^\star}(q_1,\ell):= \frac{U^{n_q}_{n_{q_1}}(n^\star,\ell)}{V^{n_q}(\ell)}\Delta Y^{n_q}(\ell); \quad q_1\neq q,$$

$$\lambda^q_{n^\star}(q_1,\ell):= \frac{U^{n_q}_{n_{q_1}}(n^\star,\ell)}{V^{n_{q_1}}(\ell)}\Delta Y^{n_{q_1}}(\ell); \quad q_1\neq q.$$

\noindent for $1\le \ell \le i$. With this notation at hand, we have
$$\langle W_q\rangle(i)=
$$

$$\lim_{n^\star\rightarrow \infty} \sum_{\ell=1}^i\sum_{q_1\neq q}\mathbb{E}[|\pi^q_{n^\star}(q_1,\ell)-\lambda^q_{n^\star}(q_1,\ell)|^2|\mathcal{F}_{\ell-1}]+$$

$$2\lim_{n^\star\rightarrow \infty}\sum_{\ell=1}^i\sum_{(q_1,q_2) \in A_q}\mathbb{E}\big[\big(\pi^q_{n^\star}(q_1,\ell)-\lambda^q_{n^\star}(q_1,\ell)\big)
\big(\pi^q_{n^\star}(q_2,\ell)-\lambda^q_{n^\star}(q_2,\ell)\big)|\mathcal{F}_{\ell-1}\big]
$$

$$=:\lim_{n^\star\rightarrow \infty} \sum_{\ell=1}^i\sum_{q_1\neq q} T_1(q_1,\ell,n^\star) + 2\lim_{n^\star\rightarrow \infty}\sum_{\ell=1}^i\sum_{(q_1,q_2) \in A_q} T_2(q_1,q_2,\ell,n^\star).$$

The martingale property and the independence of the random sample yield

\begin{equation}\label{q1}
\mathbb{E}[\Delta Y^{n_q}(\ell)\Delta Y^{n_{q_1}}(\ell)|\mathcal{F}_{\ell-1}] = \mathbb{E}[\Delta Y^{n_q}(\ell)|\mathcal{F}^{n_q}_{\ell-1}]\times\mathbb{E}[\Delta Y^{n_{q_1}}(\ell)|\mathcal{F}^{n_{q_1}}_{\ell-1}]=0,
\end{equation}

\noindent for every $q_1\neq q$. Moreover,

\begin{eqnarray}
\nonumber \mathbb{E}[|\Delta Y^{n_{q_1}}(\ell)|^2|\mathcal{F}_{\ell-1}]&=&\mathbb{E}[V^{n_{q_1}}(\ell)h^{q_1}(\ell)(1-h^{q_1}(\ell))
|\mathcal{F}^{n_{q_1}}_{\ell-1}]\\
\nonumber & &\\
\label{q2}&=& V^{n_{q_1}}(\ell)h^{q_1}[1-h^{q_1}(\ell)] \quad \forall q_1\in \mathcal{J}~\text{and}~1\le \ell\le 1.
\end{eqnarray}

From equations~(\ref{q1}) and~(\ref{q2}), we may write

\begin{eqnarray}
\nonumber T_1(q_1,\ell,n^\star)&=&\frac{|U^{n_q}_{n_{q_1}}(n^\star,\ell)|^2}{V^{n_{q_1}}(\ell)}h^{q_1}(\ell)[1-h^{q_1}(\ell)]\\
\nonumber & &\\
\label{q3}&+& \frac{|U^{n_q}_{n_{q_1}}(n^\star,\ell)|^2}{V^{n_{q}}(\ell)}h^{q}(\ell)[1-h^{q}(\ell)].
\end{eqnarray}

The same arguments for $T_1(q_1,\ell,n^\star)$ may be applied for the second term and in this case the crossing terms $T_2(q_1,q_2,\ell,n^\star)$ may be written as

\begin{eqnarray}
\nonumber T_2(q_1,q_2,\ell,n^\star) &=& \mathbb{E}[\pi^q_{n^\star}(q_1,\ell)\pi^q_{n^\star}(q_2,\ell)|\mathcal{F}_{\ell-1}]\\
\nonumber & &\\
\label{q4}&=& \frac{U^{n_q}_{n_{q_1}}(n^\star,\ell)U^{n_q}_{n_{q_2}}(n^\star,\ell)}{V^{n_q}(\ell)}h^q(\ell)[1-h^q(\ell)].
\end{eqnarray}

Summing up equations~(\ref{q3}) and~(\ref{q4}) and using assumptions \textbf{(M3)} and \textbf{(M4)}, we do have

\begin{eqnarray*}
\langle W_q \rangle(i) &=& \sum_{\ell=1}^i\sum_{q_1\neq q} \left\{\alpha^q_{q,q_1}(\ell)h^{q}(\ell)[1-h^{q}(\ell)] + \alpha^{q_1}_{q,q_1}(\ell)h^{q_1}(\ell)[1-h^{q_1}(\ell)] \right\}\\
& &\\
&+& 2\sum_{\ell=1}^i\sum_{(q_1,q_2) \in A_q}\beta^q_{q_1,q_2}(\ell)h^q(\ell)[1-h^q(\ell)].
\end{eqnarray*}
This shows that $Z_q$ has independent increments. Lemma~\ref{con_ell_2} allows us to conclude the proof.
\end{proof}

\subsection{Proof of Proposition~\ref{multlema}}

\begin{proof}
Let us fix $1\le \ell \le i$ and a weighted process $U$ which satisfies assumptions \textbf{(M1-M4)} and~\textbf{(H1)}. In view of a Cramer-Wold argument, we fix an arbitrary $a=(a_1,\ldots, a_J)\in \mathbb{R}^J$ and to shorten notation we write $\xi^{n^\star}_m(\ell):=\sum_{q=1}^J\xi^{n^\star}_{m,q}(\ell)a_q$. From Lemma~\ref{multimardiff}, we know that $\xi_{\cdot}^{n^\star}(\ell)$ is a $\mathcal{G}^{n^\star}(\ell)$-martingale difference for each $n^\star$.  Let us now check the conditions for the Central Limit theorem under martingale dependence. For any $\delta>0$, there exists a constant $C$ which only depends on $J$ and $\delta$ such that

\begin{eqnarray}
\nonumber \sum_{m=1}^{N_J} \mathbb{E}[|\xi_m^{n^\star}(\ell)|^{2+\delta} | \mathcal{G}^{n^\star}_{m-1}(\ell)]&\le& C\sum_{k=1}^J|a_k|^{2+\delta} \sum_{m=1}^{N_J}\mathbb{E}[|\xi_{m,k}^{n^\star}(i)|^{2+\delta}|\mathcal{G}^{n^\star}_{m-1}(\ell)]\\
\nonumber& &\\
\nonumber&=& C\sum_{k=1}^J|a_k|^{2+\delta} \sum_{m=1}^{N_J}\mathbb{E}[|\xi_{m,k}^{n^\star}(i)|^{2+\delta}|\mathcal{G}^{n^\star}_{m-1}(\ell)]\\
\nonumber& &\\
\label{c1}& \rightarrow& 0~\text{in probability as}~n^\star\rightarrow \infty.
\end{eqnarray}

\noindent The convergence~(\ref{c1}) is due to assumptions \textbf{(M1)} and \textbf{(M2)}. Now let us consider the predictable quadratic variation of $\sum_{m=1}^{\cdot}\xi_{m}^{n^\star}(\ell)$ at the point $N_J$ as follows

\begin{eqnarray*}
\sum_{m=1}^{N_J} \mathbb{E}\big[|\xi_{m}^{n^\star}(\ell)|^2~|\mathcal{G}^{n^\star}_{m-1}(\ell)\big]&=& \sum_{k=1}^J|a_k|^2\sum_{m=1}^{N_J}\mathbb{E}\big[|\xi_{m,k}^{n^\star}(\ell)|^2~|\mathcal{G}_{m-1}^{n^\star}(\ell)\big]\\
& &\\
&+& 2 \sum_{1\le r<k\le J}\sum_{m=1}^{N_J}a_ra_k \mathbb{E}\big[\xi_{m,k}^{n^\star}(\ell)\xi_{m,r}^{n^\star}(\ell)|\mathcal{G}^{n^\star}_{m-1}(\ell)\big]\\
& &\\
&=:& T_1(\ell,n^\star) + T_2(\ell,n^\star).
\end{eqnarray*}

Assumptions \textbf{(M3-M4)} and step~(\ref{d}) in Lemma~\ref{con_ell_2} yield

\begin{equation}\label{c2}
T_1(\ell,n^\star)\rightarrow \sum_{k=1}^J a^2_k\phi^2_k(\ell)
\end{equation}

\noindent in probability as $n^\star\rightarrow \infty$. It remains to investigate $T_2(\ell,n^\star)$. Let us fix $1\le m \le N_J$, and a pair $r\neq k$ in $\mathcal{J}$. We shall use the notation introduced in~(\ref{pi}) and~(\ref{lambda}) to write

\begin{eqnarray*}
\mathbb{E}\big[\xi_{m,k}^{n^\star}(\ell)\xi_{m,r}^{n^\star}(\ell)|\mathcal{G}^{n^\star}_{m-1}(\ell)\big]&=& \sum_{q_1\neq k}\sum_{q_2\neq r}\Bigg\{\mathbb{E}[\pi_{m,n^\star}^k(q_1,\ell)\pi_{m,n^\star}^r(q_2,\ell)| \mathcal{G}^{n^\star}_{m-1}(\ell)]\\
& &\\
&-& \mathbb{E}[\pi_{m,n^\star}^k(q_1,\ell)\lambda_{m,n^\star}^r(q_2,\ell)| \mathcal{G}^{n^\star}_{m-1}(\ell)]\\
& &\\
&-& \mathbb{E}[\lambda_{m,n^\star}^k(q_1,\ell)\pi_{m,n^\star}^r(q_2,\ell)| \mathcal{G}^{n^\star}_{m-1}(\ell)]\\
& &\\
&+& \mathbb{E}[\lambda_{m,n^\star}^k(q_1,\ell)\lambda_{m,n^\star}^r(q_2,\ell)| \mathcal{G}^{n^\star}_{m-1}(\ell)]\Bigg\}
\end{eqnarray*}

Let us fix any $q_1\neq k$ and $q_2\neq r$. Assumption \textbf{(M1)}, the independence of the random sample and~(\ref{key1}) yield

\begin{equation}\label{c3}
\mathbb{E}[\pi_{m,n^\star}^k(q_1,\ell)\pi_{m,n^\star}^r(q_2,\ell)| \mathcal{G}^{n^\star}_{m-1}(\ell)]=0.
\end{equation}

Again we may invoke assumption \textbf{(M1)}, the independence of the random sample and the definition of $\mathcal{G}^{n^\star}(\ell)$ to write the following relations

\begin{eqnarray}
\nonumber\mathbb{E}[\lambda_{m,n^\star}^k(q_1,\ell)\lambda_{m,n^\star}^r(q_2,\ell)| \mathcal{G}^{n^\star}_{m-1}(\ell)]&=&0~\text{if}~q_1\neq q_2,\\
\nonumber& &\\
\nonumber&=&V^{q_1}_m(\ell)[h^{q_1}(\ell)(1-h^{q_1}(\ell))]\\
\nonumber& &\\
\label{c4}&\times& \frac{U^{n_k}_{n_{q_1}}(n^\star,\ell)U^{n_r}_{n_{q_1}}(n^\star,\ell)}{|V^{n_{q_1}}(\ell)|^2}~\text{if}~q_1=q_2;
\end{eqnarray}

\

\begin{eqnarray}
\nonumber\mathbb{E}[\pi_{m,n^\star}^k(q_1,\ell)\lambda_{m,n^\star}^r(q_2,\ell)| \mathcal{G}^{n^\star}_{m-1}(\ell)]&=&0~\text{if}~k\neq q_2,\\
\nonumber& &\\
\nonumber&=& V^{k}_m(\ell)[h^{k}(\ell)(1-h^{k}(\ell))]\\
\nonumber& &\\
\label{c5}&\times& \frac{U^{n_k}_{n_{q_1}}(n^\star,\ell)U^{n_r}_{n_{q_2}}(n^\star,\ell)}{|V^{n_{q_2}}(\ell)|^2}~\text{if}~k=q_2;
\end{eqnarray}

\

\begin{eqnarray}
\nonumber\mathbb{E}[\lambda_{m,n^\star}^k(q_1,\ell)\pi_{m,n^\star}^r(q_2,\ell)| \mathcal{G}^{n^\star}_{m-1}(\ell)]&=&0~\text{if}~q_1\neq r,\\
\nonumber& &\\
\nonumber&=& V^{q_1}_m(\ell)[h^{r}(\ell)(1-h^{r}(\ell))]\\
\nonumber& &\\
\label{c6}&\times& \frac{U^{n_k}_{n_{q_1}}(n^\star,\ell)U^{n_r}_{n_{q_2}}(n^\star,\ell)}{|V^{n_{q_1}}(\ell)|^2}~\text{if}~q_1=r.
\end{eqnarray}

Summing up relations~(\ref{c3}),~(\ref{c4}),~(\ref{c5}) and~(\ref{c6}), we actually have

\begin{equation}\label{ast3}
T_2(\ell,n^\star) = 2\sum_{1\le r<k\le J} a_ra_k\sum_{q_1\neq k}\frac{U^{n_k}_{n_{q_1}}(n^\star,\ell) U^{n_r}_{n_{q_1}}(n^\star,\ell)}{V^{n_{q_1}}(\ell)}h^{q_1}(\ell)[1-h^{q_1}(\ell)]
\end{equation}

$$-2\sum_{1\le r< k\le J}a_ra_k\sum_{q_1\neq k}\frac{U^{n_k}_{n_{q_1}}(n^\star,\ell) U^{n_r}_{n_{k}}(n^\star,\ell)}{V^{n_{k}}(\ell)}h^{k}(\ell)[1-h^{k}(\ell)]$$

$$-2\sum_{1\le r< k\le J}a_ra_k\sum_{q_2\neq r}\frac{U^{n_k}_{n_{r}}(n^\star,\ell) U^{n_r}_{n_{q_2}}(n^\star,\ell)}{V^{n_{r}}(\ell)}h^{r}(\ell)[1-h^{r}(\ell)].$$

By making use of the assumption \textbf{(H1)}, it follows that

\begin{equation}\label{non2}
\lim_{n^\star\rightarrow \infty}T_2(\ell, n^\star)=2\sum_{1\le r < k\le J} a_ra_k \psi(k,r,\ell)
\end{equation}

\noindent in probability. Therefore, an application of the central limit theorem yields

$$\sum_{m=1}^{N_J}\xi^{n^\star}_m(\ell)\rightarrow N\Big(0, \sum_{k=1}^J a^2_k\phi^2_k(\ell) + 2\sum_{1\le r < k\le J} a_ra_k \psi(k,r,\ell)\Big),$$

\noindent weakly as $n^\star\rightarrow \infty$ and hence Lemma~\ref{mainres1} yields

$$\xi^{n^\star}(\ell)\rightarrow \tilde{Z}(\ell)~\text{weakly as}~n^\star\rightarrow \infty$$
where $\tilde{Z}(\ell) = (Z_1(\ell), \ldots, Z_J(\ell))$ (see~(\ref{zq})) has the Gaussian law $N(0,Q(\ell))$ for $1\le \ell\le i$. Relations~(\ref{c2}),~(\ref{ast1}),~(\ref{ast2}),~(\ref{ast3}) and Remark~\ref{hestimator} allow us to conclude that

$$vec\Big(\hat{Q}(n^\star, \ell)\Big)\rightarrow vec\Big(Q(\ell)\Big)$$

\noindent in probability as $n^\star\rightarrow \infty$ for $1\le \ell \le i$. It remains to check~(\ref{indvec}) but for this, we may apply the same arguments of Lemma~\ref{mainres1}. By using the notation introduced in~(\ref{wq}), let us consider

$$\tilde{W}(i)=(W_1(i),\ldots, W_J(i)).$$

By repeating the same arguments as in the proof of Lemma~\ref{mainres1}, it is straightforward to check that $\tilde{W}$ is an $\mathbb{F}$-vector martingale with independent increments on the subset $\{1,\ldots, i\}$. That is, $\tilde{Z}(j)$ and $\tilde{Z}(m)$ are independent $\mathbb{R}^J$-valued random variables for every $m\neq j$ in $\{1,\ldots, i \}$. Under these conditions we may conclude convergence~(\ref{indvec}).
\end{proof}

\section*{Acknowledgements}
This paper was completed when the second author was visiting ETH Zurich. He would like to thank the Mathematics department of ETH Zurich and Forschungsinstitut f\"{u}r Mathematik (FIM) for the very kind hospitality.



\end{document}